\documentclass[8 pt]{amsproc}
\oddsidemargin=-1cm\evensidemargin=-1cm
\begin{document}
\numberwithin{equation}{section}

\def\1#1{\overline{#1}}
\def\2#1{\widetilde{#1}}
\def\3#1{\widehat{#1}}
\def\4#1{\mathbb{#1}}
\def\5#1{\frak{#1}}
\def\6#1{{\mathcal{#1}}}

\def\C{{\4C}}
\def\R{{\4R}}
\def\N{{\4N}}
\def\Z{{\4Z}}

\author{Valentin Burcea}
\title{A normal form for a real $2$-codimensional submanifold  in
$\mathbb{C}^{N+1}$ near a  CR singularity}
\begin{abstract}
 We construct a  formal normal form for a class of real
$2$-codimensional  submanifolds $M\subset\mathbb{C}^{N+1}$ defined near a
CR singularity  approximating  the sphere. Our result gives a
 generalization of  Huang-Yin's normal form in
$\mathbb{C}^{2}$ to a higher  dimensional analogue case.
\end{abstract}

\address{School of Mathematics, Trinity College Dublin, Dublin $2$, Ireland}
\email{valentin@maths.tcd.ie}
\thanks{This
project was supported by Science Foundation Ireland grant
$06/RFP/MAT018$.}
\maketitle

\def\Label#1{\label{#1}{\bf (#1)}~}

\def\cn{{\C^n}}
\def\cnn{{\C^{n'}}}
\def\ocn{\2{\C^n}}
\def\ocnn{\2{\C^{n'}}}


\def\dist{{\rm dist}}
\def\const{{\rm const}}
\def\rk{{\rm rank\,}}
\def\id{{\sf id}}
\def\tr{{\bf tr\,}}
\def\aut{{\sf aut}}
\def\Aut{{\sf Aut}}
\def\CR{{\rm CR}}
\def\GL{{\sf GL}}
\def\Re{{\sf Re}\,}
\def\Im{{\sf Im}\,}
\def\span{\text{\rm span}}
\def\Diff{{\sf Diff}}

\def\codim{{\rm codim}}
\def\crd{\dim_{{\rm CR}}}
\def\crc{{\rm codim_{CR}}}

\def\phi{\varphi}
\def\eps{\varepsilon}
\def\d{\partial}
\def\a{\alpha}
\def\b{\beta}
\def\g{\gamma}
\def\G{\Gamma}
\def\D{\Delta}
\def\Om{\Omega}
\def\k{\kappa}
\def\l{\lambda}
\def\L{\Lambda}
\def\z{{\bar z}}
\def\w{{\bar w}}
\def\Z{{\1Z}}
\def\t{\tau}
\def\th{\theta}

\emergencystretch15pt \frenchspacing

\newtheorem{Thm}{Theorem}[section]
\newtheorem{Cor}[Thm]{Corollary}
\newtheorem{Pro}[Thm]{Proposition}
\newtheorem{Lem}[Thm]{Lemma}

\theoremstyle{definition}\newtheorem{Def}[Thm]{Definition}

\theoremstyle{remark}
\newtheorem{Rem}[Thm]{Remark}
\newtheorem{Exa}[Thm]{Example}
\newtheorem{Exs}[Thm]{Examples}

\def\bl{\begin{Lem}}
\def\el{\end{Lem}}
\def\bp{\begin{Pro}}
\def\ep{\end{Pro}}
\def\bt{\begin{Thm}}
\def\et{\end{Thm}}
\def\bc{\begin{Cor}}
\def\ec{\end{Cor}}
\def\bd{\begin{Def}}
\def\ed{\end{Def}}
\def\br{\begin{Rem}}
\def\er{\end{Rem}}
\def\be{\begin{Exa}}
\def\ee{\end{Exa}}
\def\bpf{\begin{proof}}
\def\epf{\end{proof}}
\def\ben{\begin{enumerate}}
\def\een{\end{enumerate}}
\def\beq{\begin{equation}}
\def\eeq{\end{equation}}

\section{Introduction and  Main Result}

 The study of  real submanifolds in a complex space
near a CR singularity  goes back  to the celebrated paper of
Bishop~\cite {B}. A point $p\in M$ is called a CR singularity if
it is a jumping discontinuity point for the map $M\ni q\mapsto
\dim_{\mathbb{C}}T^{0,1}_{q}M$ defined near $p$. Here
$M\subset\mathbb{C}^{N+1}$ is a real submanifold and
$T_{q}^{0,1}M$ is the CR tangent space to $M$ at $q$.

 Bishop considered the case when there exists coordinates  $(z,w)$  in
$\mathbb{C}^{2}$ such that near a CR singularity $p=0$, a real
$2$-codimensional submanifold  $M\subset\mathbb{C}^{2}$ is defined
locally by
\begin{equation}
 w=z\overline{z}+\lambda\left(z^{2}+\overline{z}^{2}\right)+\rm{O}(3),\quad\mbox{ or }\hspace{0.1
 cm}
 w=z^{2}+\overline{z}^{2}+\mbox{O}(3),
 \end{equation}
where  $\lambda\in\left[0,\infty\right]$ is a holomorphic
invariant called the Bishop invariant. When $\lambda=\infty$, $M$
is understood to be defined by the second equation from $(1)$. If
$\lambda$ is non-exceptional Moser-Webster~\cite{MW} proved
that there exists a formal transformation that sends $M$ into the
normal form
\begin{equation}
w=z\overline{z}+\left(\lambda+\epsilon
u^{q}\right)\left(z^{2}+\overline{z}^{2}\right),\quad
\epsilon\in\left\{0,-1,+1\right\},\quad q\in\mathbb{N},
\end{equation}
where $w=u+iv$.  When $\lambda=0$  Moser~\cite{Mos} constructed   the
following partial normal form:
\begin{equation}\label{moser}
 w=z\overline{z}+2\Re\left\{\displaystyle\sum_{j\geq s}a_{j}z^{j}\right\}.
\end{equation}

Here $s:=\min\left\{j\in\mathbb{N}^{\star};\hspace{0.1
cm}a_{j}\neq 0\right\}$ is the simplest higher order invariant,
known as  Moser's invariant.  Moser's partial normal form is
still subject to an infinite dimensional group action. The local
equivalence problem for the real submanifolds defined  by
(\ref{moser}) with $s<\infty$ was completely solved in a recent
 deep  paper by Huang-Yin~\cite{HY2}. Among many other
things, Huang and Yin proved that \eqref{moser} is either a
quadric or it can be formally transformed into the following
normal form defined by
\begin{equation}
w=z\overline{z}+ 2\Re\left\{\displaystyle\sum_{j\geq
s}a_{j}z^{j}\right\},\quad a_{s}=1,\quad a_{j}=0,\quad
\mbox{if}\quad j=0,1\hspace{0.1 cm}\mbox{mod}\hspace{0.1
cm}s,\quad j>s.
\end{equation}

In this paper, we construct a higher dimensional analogue of 
Huang-Yin's normal form in $\mathbb{C}^{2}$. If
$(z,w)=\left(z_{1},\dots,z_{N},w\right)$ are coordinates of
$\mathbb{C}^{N+1}$ and $M\subset\mathbb{C}^{N+1}$ a real
$2$-codimensional submanifold, we consider the case when there
exists a holomorphic change of coordinates (see  ~\cite{DTZ} or
~\cite{HY1}) such that near $p=0$, $M$ is given by
\begin{equation}w=z_{1}\overline{z}_{1}+\dots+z_{N}\overline{z}_{N}+\displaystyle\sum_{m+n\geq 3}\varphi_{m,n}(z,\overline{z}),
\end{equation}
where $\varphi_{m,n}(z,\overline{z})$ is a bihomogeneous
polynomial of bidegree $(m,n)$ in $(z,\overline{z})$.

Some of our methods extend  those from  of Huang-Yin's paper
~\cite{HY2}. First, we give a generalization of  Moser's partial
normal form \eqref{moser}, called here  Extended Moser Lemma
(Theorem~\ref{moser-lemma}), which uses the trace operator (see
e.g. ~\cite{Sh},~\cite{DZ1},~\cite{DZ2}):
\begin{equation}\mbox{tr}:=\displaystyle\sum_{k=1}^{N}\frac{\partial^{2}}{\partial
z_{k} \partial \overline{z}_{k}}.\end{equation}

In $\mathbb{C}^{2}$  Moser's partial normal form eliminates the
 terms in the local defining equation of $M$ of positive degree in both $z$ and $\overline{z}$.  The higher dimensional
case considered here brings new difficulties. In
$\mathbb{C}^{N+1}$  Extended Moser Lemma eliminates only
iterated traces of the corresponding terms. However, these terms
can still contribute to higher order terms in the construction of
the normal form. Recently, similar normal forms were constructed
for Levi-nondegenerate hypersurfaces in $\mathbb{C}^{N+1}$ by
Zaitsev in ~\cite{DZ1}. The main instrument  is given by the Fischer
decomposition.

The condition that \eqref{moser} contains nontrivial higher order
terms has the following natural generalization to the higher
dimensional case:
\begin{equation}\Re\left\{\displaystyle\sum_{k\geq 3}\varphi_{k,0}(z)\right\}\not\equiv 0,\end{equation}
where here and throughout the paper we use the following
abbreviation
$$\phi_{k,0}(z):=\phi_{k,0}(z,\bar z)$$
as the latter polynomials do not depend on $\bar z$. As a
consequence we obtain that
$s:=\mbox{min}\left\{k\in\mathbb{N}^{\star};\hspace{0.1
cm}\varphi_{k,0}(z)\not\equiv0\right\}<\infty$. Then $s$ is a
biholomorphic invariant and $\varphi_{s,0}(z)$ is invariant (as
tensor). We call the integer $s\geq3$ the generalized Moser
invariant. In this paper we will use the following notation
 \begin{equation}
 \Delta(z):=\varphi_{s,0}(z),\quad\Delta_{k}(z):=\partial_{z_{k}}\left(\varphi_{s,0}(z)\right),\quad
 k=1,\dots,N.\label{2.4}\end{equation}

 Extended Moser Lemma allows us to find just a partial normal
form. The partial normal form constructed is not unique, but that
is determined up to an action of an infinite dimensional group
${\sf Aut}_{0}\left(M_{\infty}\right)$, the formal automorphism
group
 of the quadric model
$M_{\infty}:=\left\{w=z_{1}\overline{z}_{1}+\dots+z_{N}\overline{z}_{N}\right\}$
preserving the origin. Then, the next  step is to reduce the
action by the above mentioned infinite dimensional group on the
partial normal form. In order to do this,  we use the methods
recently developed by Huang-Yin~\cite{HY2}. In particular, we
follow the ideas  of Huang-Yin's proof from ~\cite{HY2} and we 
consider  instead of the quadric model
$M_{\infty}:=\left\{w=z_{1}\overline{z}_{1}+\dots+z_{N}\overline{z}_{N}\right\}$
the model
$w=z_{1}\overline{z}_{1}+\dots+z_{N}\overline{z}_{N}+\Delta(z)+\overline{\Delta(z)}$
and by using the powerful Huang-Yin's weights system.

Before we will give the  statement of our main theorem, we
introduce the following definitions

\bd For a given homogeneous polynomial
 $V(z)=\displaystyle\sum_{\left|I\right|=k}b_{I}z^{I}$
we consider the associated Fischer differential operator
\begin{equation}V^{\star}=\displaystyle\sum_{\left|I\right|=k}\overline{b}_{I}\frac{\partial^{\left|I\right|}
}{\partial z^{I}}.\end{equation} \ed

We would like to mention that the Fischer decomposition was used
also by Ebenfelt in ~\cite{E2}.

In this paper, we consider the class of submanifolds  such that in
their defining equations, the polynomial $\Delta(z)$ defined in
(\ref{2.4})  satisfies the following nondegeneracy condition:

 \bd The polynomial $\Delta(z)$ is called nondegenerate if
for any linear forms
$\mathcal{L}_{1}(z),\dots,\mathcal{L}_{N}(z)$, one
has\begin{equation}
\mathcal{L}_{1}(z)\Delta_{1}(z)+\dots+\mathcal{L}_{N}(z)\Delta_{N}(z)\equiv0
\quad\Longrightarrow\quad \mathcal{L}_{1}(z)\equiv \dots\equiv
\mathcal{L}_{N}(z)\equiv 0 .\end{equation} \ed

In the section $2$ we prove that our non-degeneracy condition is
invariant under any linear change of coordinates.

In this paper we prove the following result:
 \bt Let $M\subset\mathbb{C}^{N+1}$ be a $2$-codimensional real (formal) submanifold
given near the point $0\in M$ by the formal power series equation
\begin{equation} w=z_{1}\overline{z}_{1}+\dots+z_{N}\overline{z}_{N}+\displaystyle\sum
_{m+n\geq3}\varphi_{m,n}(z,\overline{z}),
\label{5000}\end{equation} where $\varphi_{m,n}(z,\overline{z})$
is a bihomogeneous polynomial of bidegree $(m,n)$ in
$(z,\overline{z})$ satisfying (1.7). We assume that the
homogeneous polynomial of degree $s$ defined by (1.8) is
nondegenerate. Then there exists a unique formal map
\begin{equation}\left(z',w'\right)=\left(F(z,w),G(z,w)\right)=\left(z,w\right)+\rm{O}(2),\label{2800}\end{equation}
that transforms $M$ into the following normal form:
\begin{equation}w'=\z'_{1}\overline{z'}_{1}+\dots+z'_{N}\overline{z'}_{N}+\displaystyle\sum _{m+n\geq3 \atop m,n\neq
0}\varphi'_{m,n}\left(z',\overline{z'}\right)+2\Re\left\{\displaystyle\sum
_{k\geq s}\varphi'_{k,0}\left(z'\right)\right\},\label{3000}
\end{equation}
where $\varphi'_{m,n}\left(z',\overline{z'}\right)$ is a
bihomogeneous polynomial of bidegree $(m,n)$ in
$\left(z',\overline{z'}\right)$ satisfying the following
normalization conditions
\begin{equation}\begin{cases}
\rm{tr}^{m-1}\varphi'_{m,n}\left(z',\overline{z'}\right)=0,
 &\quad
m\leq n-1,\quad m,n\ne0 ;\cr
\rm{tr}^{n}\varphi'_{m,n}\left(z',\overline{z'}\right)=0,& \quad
m\geq n,\quad m,n\ne0.
 \end{cases}\label{CN1}
\end{equation}
\begin{equation}\left\{\begin{split}&
\left(\Delta^{t}\right)^{\star}\varphi'_{T,0}(z)=0,\quad\mbox{if}\hspace{0.2
cm} T=ts+1;\hspace{0.1
cm}t\geq1,\\&\left(\Delta_{k}\Delta^{t}\right)^{\star}\left(\varphi_{T,0}'(z)\right)=0,\hspace{0.1
cm}k=1,\dots,N,\quad\mbox{if}\hspace{0.2 cm}T=(t+1)s;\hspace{0.1
cm}t\geq1.\end{split}\right.\label{cnfnp}\end{equation} \et

A few words about the construction of the normal form. We want to
find a formal biholomorphic map sending $M$ into a formal normal
form. This leads us to study an infinite system of homogeneous
equations by truncating the original equation. As in the paper \cite {HY2} of
Huang-Yin, this system is a semi-non linear system and is very
hard to solve. We have then to use the powerful Huang-Yin's strategy
and defining the weight of $z_{k}$ to be $1$ and the weight of
$\overline{z}_{k}$ to be $s-1$, for all $k=1\dots, N$. Since ${\sf
Aut}_{0}(M_{\infty})$ is infinite-dimensional, it follows that the
homogeneous linearized normalization equations (see sections $3$
and $4$) have nontrivial kernel spaces. By using the preceding
system of weights and a similar  argument as in the paper~\cite{HY2}  of
Huang-Yin , we are able to trace precisely how the
lower order terms arise in non-linear fashion: The kernel space of
degree $2t + 1$ is restricted by imposing a normalization
condition on $\varphi'_{ts+1,0}(z)$ and the kernel space of degree
$2t + 2$ by imposing normalization conditions on
$\varphi'_{ts,0}(z)$. The non-uniqueness part of the lower degree
solutions are uniquely determined in the higher order equations.

 We would like to  mention  the  pseudo-normal form constructed  by Huang-Yin in ~\cite{HY1} as a related result to our result in $\mathbb{C}^{N+1}$.
 Our normal form is a natural generalization of  Huang-Yin's normal form in
 $\mathbb{C}^{2}$.  We  observe that our
normalization conditions are invariant under the linear changes of
coordinates that preserves the  quadric model
$M_{\infty}:=\left\{w=z_{1}\overline{z}_{1}+\dots+z_{N}\overline{z}_{N}\right\}$.
 Further related studies  concerning the  real submanifolds in the
complex space near a CR singularity were done by Ahern-Gong in ~\cite{AG}, Coffman in
~\cite{Cof1},~\cite{Cof2},~\cite{Cof3},~\cite{Cof4}, Gong
~\cite{G1},~\cite{G2},~\cite{G3}, Huang in ~\cite{HY5}, Huang-Krantz in ~\cite{HY6}.
Also, the existence of CR singularities on a compact and connected
real submanifold of codimension  $2$ was used together with natural assumtions on the boundary manifold $S$
Levi-flat hypersurfaces by Dolbeault-Tomassini-Zaitsev in ~\cite{DTZ},
~\cite{DTZ1}. Further  work in this direction was done by Lebl  in ~\cite{L1},
~\cite{L2}, ~\cite{L3}.

A few words about the paper organization: In course of section $2$
we will give a generalization of the Moser partial normal form and
we will make further preparations for our normal form
construction. The normal form construction will be presented in
the course of sections $3$ and $4$. In section $5$ we will prove
the uniqueness of the formal transformation map.

$\bf{Acknowledgements.}$ This paper was written under the
supervision of Prof. Dmitri Zaitsev. I would like to thank  him
for the introduction to the topic,
 for his patience and encouragement  during the preparation
of the manuscript. I would like  to thank Prof. Hermann Render for
pointing out the Fischer decomposition generalization from
\cite{Sh}. I am grateful also to the referees for helpful comments
on the previous versions of the manuscript.

\section{Preliminaries, notations and  Extended Moser Lemma}

Let $\left(z_{1},\dots,z_{N},w\right)$ be the coordinates from
$\mathbb{C}^{N+1}$. Let $M\subset\mathbb{C}^{N+1}$ be a real
submanifold defined near  $p=0$ by
\begin{equation}
w=z_{1}\overline{z}_{1}+\dots+z_{N}\overline{z}_{N}+\displaystyle\sum
_{m+n\geq 3}\varphi_{m,n}(z,\overline{z}), \label{ecuatiew}
\end{equation}
where $\varphi _{m,n}(z,\overline{z})$ is a bihomogeneous
polynomial of bidegree $(m,n)$ in $(z,\overline{z})$, for all
$m$,$n\geq 0$.

Let $M'$ be another submanifold defined by
\begin{equation}w'=z'_{1}\overline{z'}_{1}+\dots+z'_{N}\overline{z'}_{N}+\displaystyle\sum _{m+n\geq
3}\varphi'_{m,n}\left(z',\overline{z'}\right),
\end{equation}
where $\varphi' _{m,n}\left(z',\overline{z'}\right)$ is a
bihomogeneous polynomial of bidegree $(m,n)$ in
$\left(z',\overline{z'}\right)$, for all $m$,$n\geq 0$. We define
the hermitian product
\begin{equation}\langle z, t\rangle=z_{1}\overline{t}_{1}+\dots+z_{N}\overline{t}_{N},\quad
z=\left(z_{1},\dots,z_{N}\right),\hspace{0.1
cm}t=\left(t_{1},\dots,t_{N}\right)\in\mathbb{C}^{N}.\end{equation}

Let $\left(z',w'\right)=\left(F(z,w),G(z,w)\right)$ be a formal
map which sends $M$ to $M'$ and fixes the point
$0\in\mathbb{C}^{N+1}$. Substituting this map into $(2.2)$, we
obtain
\begin{equation} G(z,w)=\left\langle
F(z,w),F(z,w)\right\rangle+\displaystyle\sum _{m+n\geq
3}\varphi'_{m,n}\left(F(z,w),\overline{F(z,w)}\right).
\end{equation}

In the course of this paper we use the following notations
\begin{equation}\varphi_{\geq k}(z,\overline{z})=\displaystyle\sum_{m+n\geq
k}\varphi_{m,n}(z,\overline{z}),\quad
\varphi_{k}(z,\overline{z})=\displaystyle\sum_{m+n=k}\varphi_{m,n}(z,\overline{z}),\quad
k\geq 3.\label{V11}\end{equation} Substituting in (2.4)
$F(z,w)=\displaystyle\sum_{m,n\geq 0}F_{m,n}(z)w^{n}$,
$G(z,w)=\displaystyle\sum_{m,n\geq 0}G_{m,n}(z)w^{n}$, where
$G_{m,n}(z)$, $F_{m,n}(z)$ are homogeneous polynomials of degree
$m$ in $z$, by using $w$ satisfying (\ref{ecuatiew}) and notations
(\ref{V11}),   it follows that
\begin{equation}\left.\begin{split}\displaystyle &\sum
_{m,n\geq0}G_{m,n}(z)\left(\langle z,z\rangle+\varphi_{\geq
3}\right)^{n}=\left|\left|\displaystyle\sum _{m_{1},n_{1} \geq 0}
F_{m_{1},n_{1}}(z)\left(\langle z,z\rangle+\varphi_{\geq
3}\right)^{n_{1}} \right|\right|^{2}\\&+\varphi'_{\geq 3}
\left(\displaystyle\sum _{m_{2},n_{2} \geq
0}F_{m_{2},n_{2}}(z)\left(\langle z,z\rangle+\varphi_{\geq
3}\right)^{n_{2}},\overline{\displaystyle\sum _{m_{3},n_{3} \geq
0}F_{m_{3},n_{3}}(z)\left(\langle z,z\rangle+\varphi_{\geq
3}\right)^{n_{3}}}\right).
\end{split}\right.
\label{ecuatiegenerala2}\end{equation} Since our map fixes  the
point $0\in\mathbb{C}^{N+1}$, it follows that $G_{0,0}(z)=0$,
 $F_{0,0}(z)=0$. Collecting the terms of bidegree $(1,0)$  in
$(z,\overline{z})$ from (\ref{ecuatiegenerala2}), we obtain
$G_{1,0}(z)=0$.  Collecting the terms of bidegree $(1,1)$ in
$(z,\overline{z})$ from (\ref{ecuatiegenerala2}), we obtain
  \begin{equation}G_{0,1}\langle z,z\rangle=\left\langle
F_{1,0}(z),F_{1,0}(z)\right\rangle.\label{2456}\end{equation}
 Then (\ref{2456})  describes all the possible  values of
$G_{0,1}(z),F_{1,0}(z)$.  Therefore $\Im G_{0,1}=0$. By composing
with an linear automorphism of  $\Re w=\langle z,z\rangle$, we can
assume that $G_{0,1}(z)=1$, $F_{1,0}(z)=z$.

By using the same approach as in ~\cite{DZ1} (this idea was
suggested  me by Dmitri Zaitsev), the ''good'' terms that can
help us to find the formal change of coordinates under some
normalization conditions are
\begin{equation}
\varphi_{m,n}(z,\overline{z}),\quad\varphi_{m,n}'\left(z,\overline{z}\right),\quad
G_{m,n}(z)\langle z,z \rangle^{n},\quad\left\langle
F_{m,n}(z),z\right\rangle\langle z,z\rangle^{n},\quad\left\langle
z,F_{m,n}(z)\right\rangle\langle z,z\rangle^{n}.
\label{B1}\end{equation} We recall the trace decomposition (see
e.g. ~\cite{DZ1} ,~\cite{Sh}):

\bl For every bihomogeneous polynomial $P(z,\overline{z})$ and
$n\in\mathbb{N}$  there exist $Q(z,\overline{z})$ and
$R(z,\overline{z})$ unique polynomials
 such that
\begin{equation} P(z,\overline{z})=Q(z,\overline{z})\langle
z,z\rangle^{n}+R(z,\overline{z}),\quad
\rm{tr}^{n}R=0.\end{equation} \el

By using  Lemma $2.1$ and the ''good'' terms defined previously
(see (\ref{B1})) we develop a partial normal form that generalizes
 Moser's Lemma (see ~\cite{Mos}). We prove the following
statement:

\bt[Extended Moser Lemma]\label{moser-lemma} Let
$M\subset\mathbb{C}^{N+1}$ be a $2$-codimensional real-formal
submanifold. Suppose that $0\in M$ is a CR singular point and the
submanifold $M$ is defined by
\begin{equation} w=\langle
z,z\rangle+\displaystyle\sum
_{m+n\geq3}\varphi_{m,n}(z,\overline{z}),
\end{equation} where $\varphi_{m,n}\left(z,\overline{z}\right)$ is
bihomogeneous polynomial of bidegree $(m,n)$ in
$(z,\overline{z})$, for all $m,n\geq 0$. Then there exists a
unique formal map
\begin{equation}\left(z',w'\right)=\left(z+\displaystyle\sum_{m+n\geq2}F_{m,n}(z)w^{n},w+\displaystyle\sum_{m+n
\geq2}G_{m,n}(z)w^{n}\right),\label{EMLmap}\end{equation} where
$F_{m,n}(z)$, $G_{m,n}(z)$ are homogeneous polynomials in $z$ of
degree $m$ with the following normalization conditions
\begin{equation}
F_{0,n+1}(z)=0,\quad F_{1,n}(z)=0,\quad \mbox{for all}\hspace{0.1
cm} n\geq 1, \label{fn}
\end{equation}
that transforms $M$ to the following partial normal form:
\begin{equation}w'=\left\langle
z',z'\right\rangle+\displaystyle\sum _{m+n\geq3 \atop m,n\neq
0}\varphi'_{m,n}\left(z',\overline{z'}\right)+2\Re\left\{\displaystyle\sum
_{k\geq 3}\varphi'_{k,0}\left(z'\right)\right\}, \label{pnf}
\end{equation}
where $\varphi'_{m,n}(z,\overline{z})$ are bihomogeneous
polynomials of bidegree $(m,n)$  in $(z,\overline{z})$, for all
$m,n\geq 0$, that satisfy the trace normalization conditions
(\ref{CN1}).
 \et
\bpf

We construct the polynomials $F_{m',n'}(z)$ with $m'+2n'=T-1$ and
$G_{m',n'}(z)$ with $m'+2n'=T$ by induction on $T=m'+2n'$. We
assume that we have constructed the polynomials $F_{k,l}(z)$ with
$k+2l<T-1$, $G_{k,l}(z)$ with $k+2l<T$.

Collecting the terms of bidegree $(m,n)$ in $ (z,\overline{z})$
with $T=m+n$ from (\ref{ecuatiegenerala2}), we obtain
\begin{equation}
 \varphi'_{m,n}(z,\overline{z})=G_{m-n,n}(z)\left\langle
      z,z\right\rangle^{n}-\left\langle F_{m-n+1,n-1}(z),z\right\rangle\langle z,z
\rangle^{n-1}-\left\langle
 z,F_{n-m+1,m-1}(z)\right\rangle\langle z,z\rangle^{m-1}+\varphi_{m,n}(z,\overline{z})+\dots,
\label{ecg3}
\end{equation}
where ''$\dots$'' represents terms which depend on the polynomials
$G_{k,l}(z)$ with $k+2l<T$, $F_{k,l}(z)$ with $k+2l<T-1$ and on
$\varphi_{k,l}(z,\overline{z})$, $\varphi'_{k,l}(z,\overline{z})$
with  $k+l<T=m+n$.

Collecting the terms of bidegree $(m,n)$ in $ (z,\overline{z})$
with $T:=m+n\geq3$ from (\ref{ecg3}), we have to study the
following cases:

$\bf{(1)Case\hspace{0.1 cm}m<n-1,\hspace{0.1 cm}m,n\geq1.}$
Collecting the terms of bidegree $(m,n)$ in $ (z,\overline{z})$
from (\ref{ecg3}) with $m<n-1$ and $m,n\geq1$, we obtain
\begin{equation}\varphi_{m,n}'(z,\overline{z})=-\left\langle z,F_{n-m+1,m-1}(z)\right\rangle \langle
z,z\rangle^{m-1}+\dots\end{equation}

We want to use the normalization condition
$\mbox{tr}^{m-1}\varphi'_{m,n}(z,\overline{z})=0$. This allows us
to find the polynomial $F_{n-m+1,m-1}(z)$. By applying  Lemma
$2.1$ to the sum of terms which appear in ''$\dots$'', we obtain
\begin{equation}\varphi_{m,n}'(z,\overline{z})=\left(-\left\langle z,F_{n-m+1,m-1}(z)\right\rangle+D_{m,n}(z,\overline{z})\right)\langle
z,z\rangle^{m-1}+P_{1}(z,\overline{z}),\end{equation}
 where $D_{m,n}(z,\overline{z})$ is a polynomial of degree $n-m+1$ in
$\overline{z}_{1},\dots,\overline{z}_{N}$ and $1$ in
$z_{1},\dots,z_{N}$ with  determined coefficients  from the
induction hypothesis and
$\mbox{tr}^{m-1}\left(P_{1}(z,\overline{z})\right)=0$. Then, by
using the normalization condition
$\mbox{tr}^{m-1}\varphi'_{m,n}(z,\overline{z})=0$, by the
uniqueness of trace decomposition we obtain that $\left\langle
z,F_{n-m+1,m-1}(z)\right\rangle =D_{m,n}(z,\overline{z})$.  It
follows that
\begin{equation}F_{k,l}(z)=\overline{\partial_{z}\left(
D_{l+1,k+l}(z,\overline{z})\right)},\quad \mbox{for
all}\hspace{0.1 cm}k>2,l\geq 0,\label{3.7}
\end{equation}
where
$\partial_{z}:=\left(\partial_{z_{1}},\dots,\partial_{z_{N}}\right)$.

 $\bf{(2)Case\hspace{0.1 cm}m>n+1,\hspace{0.1 cm}m,n\geq1.}$
 Collecting the terms of bidegree $(m,n)$ in $
(z,\overline{z})$ from (\ref{ecg3}) with $m>n+1$ and $m,n\geq1$,
we obtain
\begin{equation}
\varphi_{m,n}'(z,\overline{z})=\left(G_{m-n,n}(z)\langle
z,z\rangle-\left\langle
F_{m-n+1,n-1}(z),z\right\rangle\right)\langle z,z
\rangle^{n-1}+\dots
\end{equation}
 In order to find the
polynomial $G_{m-n,n}(z)$ we want to use the normalization
condition $\mbox{tr}^{n}\varphi'_{m,n}(z,\overline{z})=0$. By
applying Lemma $2.1$ to the sum of terms which appear in
''$\dots$'' and to $\left\langle F_{m-n+1,n-1}(z),z\right\rangle$,
we obtain
\begin{equation}\varphi_{m,n}'(z,\overline{z})=\left(G_{m-n,n}(z)-E_{m,n}(z)\right)
\langle z,z\rangle^{n}+P_{2}(z,\overline{z}),\end{equation} where
$E_{m,n}(z)$ is a holomorphic  polynomial with determined
coefficients  by the induction hypothesis and
$\mbox{tr}^{n}\left(P_{2}(z,\overline{z})\right)=0$. Then, by
using the normalization condition
$\mbox{tr}^{n}\varphi'_{m,n}(z,\overline{z})=0$, by the uniqueness
of trace decomposition we obtain that $G_{m-n,n}(z)=E_{m,n}(z)$.
It follows that
\begin{equation}G_{k,l}(z)=E_{k+l,l}(z),\quad\mbox{for all}\hspace{0.1 cm} k\geq2,\hspace{0.1 cm}l\geq0.\end{equation}

$\bf{(3)Case\hspace{0.1 cm}(n-1,n),\hspace{0.1 cm}n\geq2.}$
Collecting the terms of bidegree $(n-1,n)$ in $ (z,\overline{z})$
from (\ref{ecg3}) with  $n\geq 2$, we obtain
\begin{equation}\varphi'_{n-1,n}(z,\overline{z})=\varphi_{n-1,n}(z,\overline{z})
-\left\langle F_{0,n-1}(z),z\right\rangle\langle
z,z\rangle^{n-1}-\left\langle z,F_{2,n-2}(z)\right\rangle\langle
z,z\rangle^{n-2}+\dots
\end{equation}

In order to find $F_{2,n-2}(z)$ we want to use the normalization
condition
  $\mbox{tr}^{n-2}\varphi'_{n-1,n}(z,\overline{z})=0$. By applying the  Lemma $2.1$ to the sum of terms from ''$\dots$'', we obtain
\begin{equation}\varphi_{n-1,n}'(z,\overline{z})=-\left(\left\langle F_{0,n-1}(z),z\right\rangle\langle z,z\rangle
+\left\langle
z,F_{2,n-2}(z)\right\rangle-C_{n-1,n}(z,\overline{z})\right)\langle
z,z\rangle^{n-2}+P_{3}(z,\overline{z}),
\end{equation}
where $\mbox{tr}^{n-2}\left(P_{3}(z,\overline{z})\right)=0$ and
$C_{n-1,n}(z,\overline{z})$ is a determined polynomial of degree
$1$ in $z_{1},\dots, z_{N}$ and degree $2$ in
$\overline{z}_{1},\dots, \overline{z}_{N}$. We take
$F_{0,n-1}(z)=0$ (see (\ref{fn})). Next, by using the
normalization condition
$\mbox{tr}^{n-2}\varphi'_{n-1,n}(z,\overline{z})=0$, by the
uniqueness of trace decomposition we obtain  that $\left\langle
z,F_{2,n-2}(z)\right\rangle=C_{n-1,n}(z,\overline{z})$. It follows
that
\begin{equation}F_{2,n-2}(z)=\overline{\partial_{z}\left(
C_{n-1,n}(z,\overline{z})\right)},\label{3.13}\end{equation} where
$\partial_{z}:=\left(\partial_{z_{1}},\dots,\partial_{z_{N}}\right)$.

$\bf{(4)Case\hspace{0.1 cm}(n,n-1),\hspace{0.1 cm}n\geq2}$.
 Collecting the terms of bidegree $(n,n-1)$ in $
(z,\overline{z})$ from (\ref{ecg3}) with $n\geq2$, we obtain
\begin{equation}\varphi_{n,n-1}'(z,\overline{z})=\left(G_{1,n-1}(z)\langle z,z\rangle-
\langle F_{2,n-2}(z),z\rangle-\left\langle
z,F_{0,n-1}(z)\right\rangle \langle z,z\rangle\right)\langle
z,z\rangle^{n-2}+\varphi_{n,n-1}(z,\overline{z})+\dots
\end{equation}

In order to find $G_{1,n-1}(z)$ we want to use the normalization
condition $\mbox{tr}^{n-1}\varphi'_{n,n-1}(z,\overline{z})=0$. By
using (\ref{fn}) and by applying Lemma $2.1$ to $\left\langle
F_{2,n-2}(z),z\right\rangle$ (see (\ref{3.13})) and to the sum of
terms from ''$\dots$'', we obtain
\begin{equation}  \varphi_{n,n-1}'(z,\overline{z})=\left(G_{1,n-1}(z)-B_{n,n-1}(z)\right)\langle
z,z\rangle^{n-1}+P_{4}(z,\overline{z}),
\end{equation}
where $\mbox{tr}^{n-1}\left(P_{4}(z,\overline{z})\right)=0$ and
$B_{n,n-1}(z)$ is a determined holomorphic polynomial. By the
uniqueness of trace decomposition we obtain that
$G_{1,n-1}(z)=B_{n,n-1}(z)$, for all  $n\geq 2$.

$\bf{(5)Case\hspace{0.1 cm}(n,n),\hspace{0.1 cm}n\geq2.}$
Collecting the  terms of bidegree $(n,n)$ in $ (z,\overline{z})$
from (\ref{ecg3}) with $n\geq2$, we obtain
\begin{equation} \varphi'_{n,n}(z,\overline{z})=G_{0,n}(z)\langle z,z\rangle^{n}-\left\langle F_{1,n-1}(z),z\right\rangle \langle z,z\rangle^{n-1}
 -\left\langle
z,F_{1,n-1}(z)\right\rangle\langle
z,z\rangle^{n-1}+\varphi_{n,n}(z,\overline{z})+\dots
\end{equation}
By taking  $F_{1,n-1}(z)=0$ (see (\ref{fn})), we obtain
$\varphi_{n,n}'(z,\overline{z})=G_{0,n}(z) \langle
z,z\rangle^{n}+\dots$. In order to find $G_{0,n}(z)$ we use the
normalization condition
$\mbox{tr}^{n}\varphi'_{n,n}(z,\overline{z})=0$. By applying the
Lemma $2.1$ to the sum of terms from ''$\dots$'' we obtain that
$\varphi_{n,n}'(z,\overline{z})=\left(G_{0,n}(z)-A_{n}\right)\langle
   z,z\rangle^{n}+P_{5}(z,\overline{z})$,
where $A_{n}$ is a determined constant and
$\mbox{tr}^{n}\left(P_{5}(z,\overline{z})\right)=0$. By the
uniqueness of trace decomposition we obtain that $G_{0,n}=A_{n}$,
for all $n\geq 3$.

$\bf{(6)Case\hspace{0.1 cm}(T,0)\hspace{0.1
cm}\mbox{and}\hspace{0.1 cm}(0,T)\hspace{0.1 cm} T\geq 3.}$
Collecting the terms of bidegree $\left(T,0\right)$ and
$\left(0,T\right)$ in $ (z,\overline{z})$ from (\ref{ecg3}), we
obtain
\begin{equation} \begin{cases} G_{T,0}(z)+\varphi'_{T,0}(z)=\varphi_{T,0}(z)+a(z)
 & \cr\varphi'_{0,T}(\bar z)=\varphi_{0,T}(\bar z)+b(\bar z) & \end{cases},\end{equation}
 where $a(z)$, $b(\bar z)$ are the sums of terms that are determined by the induction hypothesis .
By using the normalization condition $\varphi'_{0,T}(\bar
z)=\overline{\varphi'_{T,0}(z)}$ we obtain that
$G_{T,0}(z)=\varphi_{T,0}(z)+a(z)-\overline{b(\bar
z)}-\overline{\varphi_{0,T}(\bar z)}$. \epf

Extended Moser Lemma leaves undetermined an infinite number of
parameters (see (\ref{fn})). They act on the higher order terms.
In order to determine them and complete our partial normal form we
will apply in the course of sections $3$ and $4$ the following two
lemmas:

\bl Let $P(z)$ be a homogeneous pure polynomial. For every
$k\in\mathbb{N}^{\star}$, there exist $Q(z)$, $R(z)$ unique
polynomials such that
\begin{equation} P(z)=Q(z)\Delta(z)^{k}+R(z),\quad \left(\Delta^{k}\right)^{\star}\left(R(z)\right)=0.\end{equation}
\el

 \bl For every homogeneous polynomial $P(z)$ of degree $(t+1)s$
there exists a  unique decomposition
 \begin{equation}P(z)=L(z)+C(z),\quad \left(\Delta_{k}\Delta^{t}\right)^{\star}\left(C(z)\right)=0,\quad k=1,\dots,N.\end{equation}
such that
$L(z)=\left(\Delta_{1}(z)A_{1}(z)+\dots+\Delta_{N}(z)A_{N}(z)\right)\Delta(z)^{t}$,
where $A_{1}(z),\dots,A_{N}(z)$ are linear forms.
 \el

Lemmas $2.3$ and $2.4$  are consequences of the Fischer
decomposition (see ~\cite{Sh}).
\br  Lemma $2.4$ is a
particular case of  generalized Fischer's decomposition. The
polynomial $L(z)$ is uniquely determined, but the linear forms
$A_{1}(z),\dots,A_{N}(z)$ are not necessarily uniquely determined.
In order to make them uniquely determined we consider a
nondegenerate polynomial $\Delta(z)$ (see (\ref{2.4}) and
Definition $1.2$). \er

The following proposition shows us the nondegeneracy condition on
$\Delta(z)$ is invariant under any linear change of coordinates:
 \bp If $\Delta(z)$ is nondegenerate and $z\longmapsto Az$ is a linear
 change of coordinates, then $\Delta(Az)$ is also nondegenerate.
 \ep
 \begin{proof}

Let $\widetilde{\Delta}(z)=\Delta\left(Az\right)$, where
$A=\left\{a_{jk}\right\}_{1\leq j,k \leq N}$. Therefore
$\widetilde{\Delta}_{j}(z)=\displaystyle\sum_{k=1}^{N}\Delta_{k}\left(Az\right)a_{jk}$,
for all $j=1,\dots, N$. We consider
$\mathcal{L}_{1}(z),\dots,\mathcal{L}_{N}(z)$ linear forms such
that
$\mathcal{L}_{1}(z)\widetilde{\Delta}_{1}(z)+\dots+\mathcal{L}_{N}(z)\widetilde{\Delta}_{N}(z)\equiv0,$
or equivalently
$\displaystyle\sum_{j,k=1}^{N}\Delta_{k}\left(Az\right)\mathcal{L}_{j}(z)a_{jk}\equiv
0$. Since $\Delta(z)$ is nondegenerate and $\left\{
a_{jk}\right\}_{1\leq j,k \leq N}$ is invertible it follows that
$\mathcal{L}_{1}(z)\equiv\dots \equiv\mathcal{L}_{N}(z)\equiv 0$.
 \end{proof}

$\bf{The\hspace{0.1 cm}system\hspace{0.1 cm}of\hspace{0.1
cm}weights:}$  Following Huang-Yin's ideas  from  ~\cite{HY2}, we
define the  system of weights for
$z_{1},\overline{z}_{1},\dots,z_{N},\overline{z}_{N}$ as follows.
We define $\mbox{wt}\left\{z_{k}\right\}=1$ and
$\mbox{wt}\left\{\overline{z}_{k}\right\}=s-1$, for all
$k=1,\dots,N$. If $A(z,\overline{z})$ is a formal power series we
write $\mbox{wt}\left\{A(z,\overline{z})\right\}\geq k$ if
$A\left(tz,t^{s-1}\overline{z}\right)=\mbox{O}\left(t^{k}\right)$.
We also write $\mbox{Ord}\left\{A(z,\overline{z})\right\}=k$ if
$A\left(tz,t\overline{z}\right)=t^{k} A(z,\overline{z})$.  We
denote by $\Theta_{m}^{n}(z,\overline{z})$ a series in
$(z,\overline{z})$ of weight at least $m$ and order at least $n$.
In the particular case when $\Theta_{m}^{n}(z,\overline{z})$ is
just a polynomial we use the notation
$\mathbb{P}_{m}^{n}(z,\overline{z})$. We define the set of the
normal weights
$$\mbox{wt}_{nor}\left\{w\right\}=2,\quad
\mbox{wt}_{nor}\left\{z_{1}\right\}=\dots=\mbox{wt}_{nor}\left\{z_{N}\right\}
=\mbox{wt}_{nor}\left\{\overline{z}_{1}\right\}=\dots=\mbox{wt}_{nor}\left\{\overline{z}_{N}\right\}=1.$$

$\bf{Notations:}$ If $h(z,w)$ is a formal power series with no
constant term we introduce the following notations
\begin{equation}\left.\begin{split}&h(z,w)=\displaystyle\sum_{l\geq
1}h_{nor}^{(l)}(z,w),\quad\mbox{where}\hspace{0.1
cm}\mbox{$h_{nor}^{(l)}\left(tz,t^{2}w\right)=t^{l}h_{nor}^{(l)}(z,w)$},\\&
\mbox{$h_{\geq l}(z,w)=\displaystyle\sum_{k\geq
l}h_{nor}^{(k)}(z,w)$}.\end{split}\right.\label{2.11}\end{equation}

\section{Proof of Theorem $1.3$-Case $T+1=ts+1$, $t\geq1$}

By applying Extended Moser Lemma  we can assume that $M$ is given
by the following equation
\begin{equation} w=\langle
z,z\rangle+\displaystyle\sum_{m+n
\geq3}^{T+1}\varphi_{m,n}(z,\overline{z})+\mbox{O}\left(T+2\right),
\label{ec4}\end{equation} where $\varphi_{m,n}(z,\overline{z})$
satisfies (\ref{CN1}), for all $3\leq m+n\leq T$.

We perform induction on $T\geq 3$. Assume that (\ref{cnfnp}) holds
for
 $\varphi_{k,0}(z)$, for all $k=s+1,\dots,T$ with $k=0,1\hspace{0.1 cm}\mbox{mod}\hspace{0.1 cm}(s)$. If
$T+1\not\in \left\{ts;\hspace{0.1 cm}
t\in\mathbb{N}^{\star}-\left\{1,2\right\}
\right\}\cup\left\{ts+1;\hspace{0.1
cm}t\in\mathbb{N}^{\star}\right\}$ we apply Extended Moser Lemma.
In the case when $T+1\in\left\{ts;\hspace{0.1 cm
t\in\mathbb{N}^{\star}}-\{1\}\right\}\cup\left\{ts+1;\hspace{0.1
cm}t\in\mathbb{N}^{\star}\right\}$,  we will look for a formal map
which sends our submanifold $M$ to a new submanifold $M'$ given by
\begin{equation}w'=\left\langle z',z'\right\rangle+\displaystyle\sum_{m+n\geq3}^{T+1}\varphi'_{m,n}\left(z',\overline{z'}\right)+\mbox{O}\left(T+2\right)
,\label{em2}\end{equation} where
$\varphi'_{m,n}\left(z',\overline{z'}\right)$ satisfies
(\ref{CN1}), for all $3\leq m+n\leq T$ and
$\varphi'_{k,0}\left(z'\right)$ satisfies (\ref{cnfnp}), for all
$k=s+1,\dots,T$ with $k=0,1\hspace{0.1 cm}\mbox{mod}\hspace{0.1
cm}(s)$. We will obtain that $\varphi'_{k,0}(z)=\varphi_{k,0}(z)$
for all $k=s,\dots,T$.

In the course of this section we consider the case when
$T+1=ts+1$. We are looking for a biholomorphic transformation of
the following type
 \begin{equation}\left.\begin{split}&\quad\quad\quad\quad\quad\left(z',w'\right)=\left(z+F(z,w),w+G(z,w)\right)\\&F(z,w)=\displaystyle\sum_{l=0}^{T-2t}F_{nor}^{(2t+l)}(z,w),\quad
G(z,w)=\displaystyle\sum_{\tau=0}^{T-2t}G_{nor}^{(2t+1+\tau)}(z,w)\end{split}\right.,
\label{t1}
\end{equation} that maps $M$ into $M'$ up to the order $T+1=ts+1$. In order for the
preceding mapping to be uniquely determined we assume that
$F_{nor}^{(2t+l)}(z,w)$  is normalized as in Extended Moser Lemma,
for all $l=1,\dots T$. Substituting (\ref{t1}) into (\ref{em2}) we
obtain
\begin{equation}w+G(z,w)= \left\langle z+F(z,w), z+F(z,w) \right\rangle +
\displaystyle\sum_{m+n\geq3}^{T+1}\varphi'_{m,n}\left(z+F(z,w),\overline{z+F(z,w)}\right)
 +\mbox{O}\left(T+2\right),\label{4.4}\end{equation}
 where $w$ satisfies (\ref{ec4}). By making some
simplifications in (\ref{4.4}) by using (\ref{ec4}), we obtain
\begin{equation}\left.
\begin{split}  \displaystyle\sum_{\tau=0}^{T-2t}G_{nor}^{(2t+1+\tau)}\left(z,\langle z, z\rangle+\varphi_{\geq3}(z,\overline{z})\right)&= 2\Re\left\langle
z,\displaystyle\sum_{l=0}^{T-2t}F_{nor}^{(2t+l)}\left(z,\langle
z,z\rangle+\varphi_{\geq3}(z,\overline{z})\right)\right\rangle+\left|\left|\displaystyle\sum_{l=0}^{T-2t}F_{nor}^{(2t+l)}\left(z,\langle
z,
z\rangle+\varphi_{\geq3}(z,\overline{z})\right)\right|\right|^{2}\\&\quad+
\varphi'_{\geq3}\left(z+\displaystyle\sum_{l=0}^{T-2t}F_{nor}^{(2t+l)}\left(z,\langle
z,
z\rangle+\varphi_{\geq3}(z,\overline{z})\right),\overline{z+\displaystyle\sum_{l=0}^{T-2t}F_{nor}^{(2t+l)}\left(z,\langle
z, z\rangle+\varphi_{\geq3}(z,\overline{z})\right)}\right)\\&
\quad-\varphi_{\geq3}(z,\overline{z}).
\end{split}\right.\label{ecg4}
\end{equation}

Collecting the terms with the same bidegree from (\ref{ecg4}), we
find $F(z,w)$ and $G(z,w)$ by applying Extended Moser Lemma. Since
we don't have components of $F(z,w)$ of normal weight less than
$2t$ and $G(z,w)$ with normal weight less than $2t+1$, collecting
in (\ref{ecg4}) the terms with the same bidegree $(m,n)$ in
$(z,\overline{z})$ with $m+n<2t+1$, we obtain that
$\varphi'_{m,n}(z,\overline{z})=\varphi_{m,n}(z,\overline{z})$.

Collecting the terms of bidegree $(m,n)$ in $(z,\overline{z})$
with $m+n=2t+1$ (like in the
 Extended Moser Lemma's proof) we find $G_{nor}^{(2t+1)}(z,w)$
and $F_{nor}^{(2t)}(z,w)$ as follows. We make the following claim:
\bl $G_{nor}^{(2t+1)}(z,w)=0$,
$F_{nor}^{(2t)}(z,w)=aw^{t}-z\langle z,a \rangle w^{t-1}$, where
$a=\left(a_{1},\dots,a_{N}\right)\in\mathbb{C}^{N}$.
 \el
\begin{proof}

Collecting the pure terms of degree $2t+1$  from (\ref{ecg4}), we
obtain that $\varphi_{2t+1,0}(z)=\varphi'_{2t+1,0}(z)$. Collecting
the terms of bidegree $(m,n)$ with $m+n=2t+1$ in
$(z,\overline{z})$ and $0<m<n-1$ (\ref{ecg4}), we obtain
\begin{equation}\varphi_{m,n}'(z,\overline{z})=-\left\langle
z,F_{n-m+1,m-1}(z) \right\rangle \langle
z,z\rangle^{m-1}+\varphi_{m,n}(z,\overline{z}).\end{equation}
 Since $\varphi_{m,n}(z,\overline{z})$, $\varphi'_{m,n}(z,\overline{z})$
 satisfy (\ref{CN1}), by the uniqueness of the trace decomposition, we
 obtain $F_{n-m+1,m-1}(z)=0$. Collecting the terms of bidegree $(m,n)$ in
$(z,\overline{z})$ with $m+n=2t+1$ and $m>n+1$ from (\ref{ecg4}),
we obtain
\begin{equation}\varphi_{m,n}'(z,\overline{z})=G_{m-n,n}(z)\langle
z,z\rangle^{n}-\left\langle F_{m-n+1,n-1}(z),z\right\rangle
\langle z,z
\rangle^{n-1}+\varphi_{m,n}(z,\overline{z}).\end{equation} Since
$F_{m-n+1,n-1}(z)=0$ it follows that $G_{m-n,n}(z)=0$. Collecting
the  terms of bidegree $(t-1,t)$ and  $\left(t,t-1\right)$ in
$(z,\overline{z})$ from (\ref{ecg4}), we obtain the following two
equations
\begin{equation}\left.\begin{split}&\quad\quad\hspace{0.2 cm}\varphi_{t-1,t}'(z,\overline{z})=-\left(\left\langle
F_{0,t-1}(z),z\right\rangle \langle z,z\rangle+\left\langle
z,F_{2,t-2}(z)\right\rangle\right)\langle z,z\rangle^{t-2}+
\varphi_{t-1,t}(z,\overline{z}),\\&\varphi_{t,t-1}'(z,\overline{z})=G_{1,t-1}(z)\langle
z,z\rangle^{t-1}-\left(\left\langle F_{2,t-2}(z),z\right\rangle
+\left\langle z,F_{0,t-1}(z)\right\rangle\langle
z,z\rangle\right)\langle
z,z\rangle^{t-2}+\varphi_{t,t-1}(z,\overline{z}).\end{split}\right.\label{90000}\end{equation}
By using (\ref{90000}) it follows that $G_{1,t-1}(z)=0$. We set
$F_{0,t-1}(z)=a=:\left(a_{1},\dots,a_{N}\right)$ and we write
$F_{2,t-2}(z)=\left(F_{2,t-2}^{1}(z),\dots,F_{2,t-2}^{N}(z)\right)$.
Since $\varphi_{m,n}(z,\overline{z}),\hspace{0.1
cm}\varphi'_{m,n}(z,\overline{z})$ satisfy (\ref{CN1}), by the
uniqueness of the trace decomposition, from (\ref{90000}) we
obtain the equation $\langle z,a\rangle \langle z,z\rangle
+\left\langle F_{2,t-2}(z),z\right\rangle=0$, that can be solved
as
\begin{equation}F_{2,t-2}^{k}(z)=-\frac{\partial}{\partial
\overline{z}_{k}}\left(\langle z,a\rangle\langle z,z
\rangle\right)=-z_{k}\langle z,a\rangle,\quad
k=1,\dots,N.\end{equation} Therefore
$F_{nor}^{(2t)}(z,w)=aw^{t}-z\langle z,a \rangle
w^{t-1},\hspace{0.1 cm}\mbox{where}\hspace{0.1
cm}a=\left(a_{1},\dots,a_{N}\right)\in\mathbb{C}^{N}.$
\end{proof}
By Lemma $3.1$ we conclude that
$F(z,w)=F_{nor}^{(2t)}(z,w)+F_{\geq 2t+1}(z,w)$ and
$G(z,w)=G_{\geq 2t+2}(z,w)$ (see (\ref{2.11})). We also have
$F_{\geq
2t+1}(z,w)=\displaystyle\sum_{k+2l\geq2t+1}F_{k,l}(z)w^{l}$, where
$F_{k,l}(z)$ is a homogeneous polynomial of degree $k$. It follows
that \begin{equation}\mbox{wt}\left\{F_{\geq
2t+1}(z,w)\right\}\geq\displaystyle\min_{k+2l\geq2t+1}\{k+ls\}
\geq\displaystyle\min_{k+2l\geq2t+1}\{k+2l\}\geq2t+1.\end{equation}
Next, we prove that $\mbox{wt}\left\{\overline{F_{\geq
2t+1}(z,w)}\right\}\geq ts+s-1$. Since
$\mbox{wt}\left\{\overline{F_{\geq
2t+1}(z,w)}\right\}\geq\displaystyle\min_{k+2l\geq2t+1}\left\{k(s-1)+ls\right\}$,
it is enough to  prove that $k(s-1)+ls\geq ts+s-1$ for $k+2l\geq
2t+1$. Since we can write the latter inequality as
$(k-1)(s-1)+ls\geq ts$, for $(k-1)+2l\geq 2t$, it is enough to
prove that $k(s-1)+ls\geq ts$, for $k+2l\geq 2t$. Since $s\geq3$
it follows that $ks-2k\geq0$. Hence $2k(s-1)+2ls\geq ks+2ls$. It
follows that $k(s-1)+ls\geq\frac{s}{2}(k+2l)\geq
\frac{2ts}{2}=ts$.

 \bl By using  the previous
calculations, we give the following immediate estimates
 \begin{equation}\left.\begin{split}& \rm{wt}\left\{F_{\geq 2t+1}(z,w) \right\}\geq 2t+1
,\quad \rm{wt} \left\{\overline{F_{\geq 2t+1}(z,w)} \right\}\geq
ts+s-1,\quad\rm{wt}\left\{\left\|F_{\geq
2t+1}(z,w)\right\|^{2}\right\}\geq
ts+2,\\&\quad\quad\rm{wt}\left\{F_{nor}^{(2t)}(z,w)\right\}\geq
ts+2-s,\quad \rm{wt} \left\{
\overline{F_{nor}^{(2t)}(z,w)}\right\}\geq ts,\quad
\rm{wt}\left\{\left\|F_{nor}^{(2t)}(z,w)\right\|^{2}\right\}\geq
ts+2,\\&\quad\quad\quad\rm{wt}\left\{\left\langle
F_{nor}^{(2t)}(z,w),F_{\geq 2t+1}(z,w)\right\rangle
\right\},\quad\rm{wt}\left\{\left\langle F_{\geq
2t+1}(z,w),F_{nor}^{(2t)}(z,w)\right\rangle \right\}\geq ts+2,
\end{split}\right.\label{1e1}\end{equation}
where $w$ satisfies (\ref{ec4}). \el As a consequence of the
preceding estimates, we obtain
 \begin{equation}\left\|F(z,w)\right\|^{2}=\left\|F_{nor}^{(2t)}(z,w)\right\|^{2}+2\Re\left\langle
F_{nor}^{(2t)}(z,w),F_{\geq 2t+1}(z,w)\right\rangle+
\left\|F_{\geq
 2t+1}(z,w)\right\|^{2}=\Theta_{ts+2}^{2t+2}(z,\overline{z}),\label{wt1}\end{equation}
 where $w$ satisfies (\ref{ec4}). We observe that the preceding  power series
 $\Theta_{ts+2}^{2t+2}(z,\overline{z})$ has the  property $\rm{wt}\left\{\overline{\Theta_{ts+2}^{2t+2}(z,\overline{z})}\right\}\geq
 ts+2$.

 In order to apply Extended Moser Lemma in (\ref{ecg4}) we
 have to evaluate the weight and the order of the terms which appear and are not ''good''.  Beside the
previous weight estimates (see (\ref{1e1}) and (\ref{wt1})) we
also need to prove the following lemmas:

 \bl For all  $m,n\geq 1$ and $w$ satisfying (\ref{ec4}), we
have the following estimate
\begin{equation}\varphi'_{m,n}\left(z+F(z,w),\overline{z+F(z,w)}\right)=\varphi'_{m,n}(z,\overline{z})
+2\Re\left\langle \Theta_{s}^{2}(z,\overline{z}),
\overline{F_{\geq
2t+1}(z,w)}\right\rangle+\Theta_{ts+2}^{2t+2}(z,\overline{z}),\end{equation}
where
$\rm{wt}\left\{\overline{\Theta_{ts+2}^{2t+2}(z,\overline{z})}\right\}\geq
 ts+2$.
 \el
\begin{proof}
We make the  expansion
$\varphi'_{m,n}\left(z+F(z,w),\overline{z+F(z,w)}\right)
=\varphi'_{m,n}(z,\overline{z})+\dots ,$ where in ''$\dots$'' we
have different types of terms involving $F_{k',l'}(z)$ with
$k'+2l'<m+n$  and normalized terms
$\varphi_{k,l}(z,\overline{z})$, $\varphi'_{k,l}(z,\overline{z})$
with $k+l<m+n$. In order to study the weight and the order  of
terms which can appear in ''$\dots$'' it is enough to study the
weight and the order of the following particular terms
$$A_{1}(z,w)=F_{1}(z,w)z^{I}\overline{z}^{J},\quad A_{2}(z,w)=
z^{I_{1}}\overline{z}^{J_{1}}\overline{F_{1}(z,w)} ,\quad
B_{1}(z,w)=F_{2}(z,w)z^{I}\overline{z}^{J},\quad
B_{2}(z,w)=\overline{F_{2}(z,w)}z^{I_{1}}\overline{z}^{J_{1}},$$
where $F_{1}(z,w)$ is the first component of $F_{nor}^{(2t)}(z,w)$
and $F_{2}(z,w)$ is the first component of $F_{\geq 2t+1}(z,w)$.
Here we assume  that $\left|I\right|=m-1$, $\left|I_{1}\right|=m$,
$\left|J_{1}\right|=n-1$, $\left|J\right|=n$.

By using (\ref{1e1}) we obtain
$\mbox{wt}\left\{A_{1}(z,w)\right\}\geq m-1+ts+2-s+n(s-1)\geq
ts+2$. It is equivalent to prove that $m-1+s(n-1)-n\geq 0$. This
is true because $m-1+s(n-1)-n\geq m-1+3(n-1)-n\geq m+3n-4-n\geq
3+n-4\geq0$. On the other hand, we have
$\mbox{Ord}\left\{A_{1}(z,w)\right\}\geq m-1+2t+n\geq 2t+2$.

 By using (\ref{1e1}) we obtain
$\mbox{wt}\left\{A_{2}(z,w)\right\}\geq m+ts+(n-1)(s-1)\geq
ts+2\Longleftrightarrow m+(s-1)(n-1)\geq 2$. We have
$m+(n-1)(s-1)\geq m+2(n-1)\geq m+2n-4\geq 0$, and this is true
because $m+n\geq3$ and  $m,n\geq1$. On the other hand we have
$\mbox{Ord}\left\{A_{2}(z,w)\right\}\geq m+2t+n-1\geq 2t+2$.

In the same  way we obtain that
$\rm{Ord}\left\{B_{1}(z,w)\right\}$,
$\rm{Ord}\left\{B_{2}(z,w)\right\}\geq 2t+1$. By using
(\ref{1e1}), every term from ''$\dots$'' that depends on
$F_{2}(z,w)$ can be written as
$\Theta_{s}^{2}(z,\overline{z})F_{2}(z,w)$. From here we obtain
our claim.
\end{proof}

 \bl For $w$ satisfying (\ref{ec4}) and for all $k>s$, we have the following estimation
 \begin{equation}\varphi'_{k}\left(z+F(z,w)\right)=\varphi'_{k}(z)+2\Re \left\langle
\Theta_{s}^{2}(z,\overline{z}),\overline{ F_{\geq
2t+1}(z,w)}\right\rangle+\Theta_{ts+2}^{2t+2}(z,\overline{z}),\end{equation}
where
$\rm{wt}\left\{\overline{\Theta_{ts+2}^{2t+2}(z,\overline{z})}\right\}\geq
 ts+2$.
 \el
\begin{proof}

We make the expansion
$\varphi'_{k}\left(z+F(z,w)\right)=\varphi'_{k}(z)+\dots$ . In
order to study the weight and the order  of terms which can appear
in ''$\dots$'' it is enough to study the weight and the order of
the following terms
$$A(z,w)=F_{1}(z,w)z^{I},\quad B(z,w)=F_{2}(z,w)z^{I},$$ where $F_{1}(z,w)$ is the
first component of $F_{nor}^{(2t)}(z,w)$ and $F_{2}(z,w)$ is the
first component of $F_{\geq 2t+1}(z,w)$. Here we assume  that
$\left|I\right|=m-1\geq s$. Then, by (\ref{1e1}), we obtain that
$\mbox{wt}\left\{A(z,w)\right\}\geq s+ts+2-s\geq ts+2$. On the
other hand, we have $\mbox{Ord}\left\{A(z,w)\right\}\geq s+2t\geq
2t+2$. By using (\ref{1e1}), every term from ''$\dots$'' that
depends on $F_{2}(z,w)$ can be written as
$\Theta_{s}^{2}(z,\overline{z})F_{2}(z,w)$. From here we obtain
our claim.
\end{proof}

We want to evaluate the weight and the order of the other terms of
(\ref{ecg4}). By Lemma $4.3$ and Lemma $4.4$, it remains to
evaluate the order and the weight of the terms of the following
expression
\begin{equation}\left.\begin{split} S(z,\overline{z})&=2\Re\left \langle F(z,w),z\right \rangle
+2\Re\left\{\varphi'_{s}\left(z+F(z,w)\right)\right\},\\&
=2\Re\left \langle F^{(2t)}_{nor}(z,w)+F_{\geq 2t+1}(z,w),z\right
\rangle+2\Re\left\{\Delta\left(z+F^{(2t)}_{nor}(z,w)+F_{\geq
2t+1}(z,w)\right)\right\},
\end{split}\right.\label{18}\end{equation}
where $w$ satisfies (\ref{ec4}).

 \bl For $F_{nor}^{(2t)}(z,w)$
given by Lemma $4.1$ and $w$ satisfying (\ref{ec4}) we have
\begin{equation}2\Re\left \langle F_{nor}^{(2t)}\left(z,w\right),z\right
\rangle=2\Re\left\{\langle
z,a\rangle\Delta(z)w^{t-1}\right\}+\Theta_{ts+2}^{2t+2}(z,\overline{z}),\end{equation}
where
$\rm{wt}\left\{\overline{\Theta_{ts+2}^{2t+2}(z,\overline{z})}\right\}\geq
 ts+2$.
 \el
\begin{proof}

We  compute
\begin{equation}\left.\begin{split}
2\Re\left \langle F_{nor}^{(2t)}\left(z,w\right),z\right
\rangle&=2\Re \left\{w^{t}\left \langle a ,z\right\rangle\right\}-
2\Re \left\{ \langle z,a \rangle\langle z,z\rangle
 w^{t-1}\right\},\\&=2\Re \left\{\left \langle z ,a\right\rangle w^{t}- \langle z,a \rangle\langle z,z\rangle
 w^{t-1}\right\}+\langle a,z\rangle\left(w^{t}-\overline{w}^{t}\right) +\langle z,a\rangle\left(\overline{w}^{t}-w^{t}\right),\\&=2\Re\left\{\langle
z,a\rangle\Delta(z)w^{t-1}\right\}+\Theta_{ts+2}^{2t+2}(z,\overline{z}),\end{split}\right.\label{22}\end{equation}
where
$\mbox{wt}\left\{\overline{\Theta_{ts+2}^{2t+2}(z,\overline{z})}\right\}\geq
 ts+2$.
 \end{proof}
In the course of our proof we  will use the notation
$\Delta'(z)=\left(\Delta_{1}(z),\dots,\Delta_{N}(z)\right)$. It
remains to prove the following lemma

 \bl For $w$ satisfying  (\ref{ec4})
we have the following estimate
\begin{equation}\left.\begin{split}2\Re\{\Delta\left(z+F(z,w)\right)\}=&2\Re\left\{\Delta(z)-s\langle
z,a\rangle\Delta(z)w^{t-1}\right\}\\&+2\Re\left\langle
\Delta'(z)+\Theta_{s}^{2}(z,\overline{z}),\overline{F_{\geq
2t+1}(z,w)}\right\rangle+\Theta_{ts+2}^{2t+2}(z,\overline{z}),\end{split}\right.\label{24}\end{equation}
where
$\rm{wt}\left\{\overline{\Theta_{ts+2}^{2t+2}(z,\overline{z})}\right\}\geq
 ts+2$.\el
\begin{proof}

 By using the Taylor expansion it follows that
\begin{equation}\left.\begin{split}2\Re\left\{\Delta\left(z+F(z,w)\right)\right\}=2\Re\left\{\Delta(z)+
\displaystyle\sum_{k=1}^{N}\Delta_{k}(z)F_{\geq
2t}^{k}(z,w)+L(z,\overline{z})\right\},
\end{split}\right.\label{25}\end{equation}
where $F_{\geq 2t}^{k}(z,w)=\left(F^{1}_{\geq
2t}(z,w),\dots,F^{N}_{\geq 2t}(z,w)\right)$ and
$L(z,\overline{z})=\left\langle\Theta_{s}^{2}(z,\overline{z}),
\overline{F_{\geq2t+1}(z,w)}\right\rangle$.  We compute
\begin{equation}\left.\begin{split} \displaystyle\sum_{k=1}^{N}2\Re\left\{\Delta_{k}(z)F_{\geq 2t}^{k}(z,w)\right\}&
=\displaystyle
\displaystyle\sum_{k=1}^{N}2\Re\left\{\Delta_{k}(z)\left(a_{k}w^{t}-z_{k}\langle
z,a \rangle w^{t-1}+F_{\geq2t+1}^{k}(z,w)\right)\right\},
 \\&=\Theta_{ts+2}^{2t+2}(z,\overline{z})-2s\Re\left\{\langle z,a \rangle\Delta(z)
w^{t-1}\right\}+2\Re\left\langle \Delta'(z),\overline{F_{\geq
2t+1}(z,w)}\right\rangle,\end{split}\right.\label{26}\end{equation}
where
 $\rm{wt}\left\{\overline{\Theta_{ts+2}^{2t+2}(z,\overline{z})}\right\}\geq
 ts+2$.
\end{proof}

 For $w$ satisfying (\ref{ec4}), by Lemma $3.5$ and by Lemma $3.6$, we can rewrite (\ref{18}) as follows
\begin{equation}\left.\begin{split}S(z,\overline{z})=2(1-s)\Re\left\{\langle z,a \rangle\Delta(z)w^{t-1}\right\}+2\Re\left\langle
\overline{z}+\Delta'(z)+\Theta_{s}^{2}(z,\overline{z}),\overline{
F_{\geq
2t+1}\left(z,w\right)}\right\rangle+\Theta_{ts+2}^{2t+2}(z,\overline{z}),\end{split}\right.\end{equation}
where
$\rm{wt}\left\{\overline{\Theta_{ts+2}^{2t+2}(z,\overline{z})}\right\}\geq
 ts+2$. By Lemmas $3.1-3.6$ we obtain
\begin{equation}\left.
\begin{split}
G_{\geq 2t+2}\left(z,\langle
z,z\rangle+\varphi_{\geq3}(z,\overline{z})\right)=&2(1-s)\Re\left\{\langle
z,a \rangle\Delta(z)\left(\langle z,z\rangle
+\varphi_{\geq3}(z,\overline{z})\right)^{t-1}\right\}
\\&+2\Re\left\langle \overline{z}+\Delta'(z)+\Theta_{s}^{2}(z,\overline{z}),\overline{F_{\geq
2t+1}\left(z,\langle
z,z\rangle+\varphi_{\geq3}(z,\overline{z})\right)}\right\rangle\\&+\varphi_{\geq2t+2}(z,\overline{z})-\varphi'_{\geq2t+2}(z,\overline{z})+\Theta_{ts+2}^{2t+2}(z,\overline{z}),
\end{split}\right.\label{27}\end{equation}
where
$\rm{wt}\left\{\overline{\Theta_{ts+2}^{2t+2}(z,\overline{z})}\right\}\geq
 ts+2$.

 Assume that $t=1$. Collecting the terms of total degree $k<s+1$ in $(z,\overline{z})$  from (\ref{27}) we
 find the polynomials
$\left(G_{nor}^{(k+1)}(z,w),F_{nor}^{(k)}(z,w)\right)$ for all
$k<s$. Collecting the terms of total degree $m+n=s+1$ in
$(z,\overline{z})$ from (\ref{27}), we  obtain
\begin{equation}\left.
\begin{split}G_{nor}^{(s+1)}\left(z,\langle z,z \rangle\right)=&2(1-s)\Re\left\{\langle z,a
\rangle\Delta(z)\right\} +2\Re\left\langle
z,F_{nor}^{(s)}\left(z,\langle z,z\rangle\right)\right
\rangle+\varphi'_{s+1}(z,\overline{z})-\varphi_{s+1}(z,\overline{z})+\left(\Theta_{1}\right)_{s+2}^{s+1}(z,\overline{z}).
\end{split}\right.\label{28}
\end{equation} By applying Extended Moser Lemma  we find a solution
$\left(G_{nor}^{(s+1)}(z,w),F_{nor}^{(s)}(z,w)\right)$ for the
latter equation. We consider the following  Fischer decompositions
\begin{equation}\varphi_{s+1,0}(z)=Q(z)\Delta(z)+R(z),\quad
\varphi'_{s+1,0}(z)=Q'(z)\Delta(z)+R'(z),\label{29}\end{equation}
where
$\Delta^{\star}\left(R(z)\right)=\Delta^{\star}\left(R'(z)\right)=0$.
We want to put the normalization condition
 $\Delta^{\star}\left(\varphi'_{s+1,0}(z)\right)=0$. Collecting the pure terms of degree $s+1$ in (\ref{28}), by (\ref{29}) we
 obtain
 \begin{equation}\varphi'_{s+1,0}(z)=\varphi_{s+1,0}(z)-(1-s)\langle z,a \rangle \Delta(z)=\left(Q(z)-(1-s)\langle
 z,a\rangle\right)\Delta(z)+R(z),\label{f01}\end{equation}
where $Q(z)$ is a determined polynomial  of degree $1$ in
$z_{1},\dots,z_{N}$. It follows that $Q'(z)= Q(z)-(1-s)\langle
z,a\rangle$ and $R'(z)=R(z)$. Then the normalization condition
 $\Delta^{\star}\left(\varphi'_{s+1,0}(z)\right)=0$ is equivalent to finding $a$ such that
$Q'(z)= Q(z)-(1-s)\langle z,a\rangle=0.$ The last equation
provides us the free parameter $a$.

 Assuming that $t\geq 2$, we prove the following lemma (this is the analogue of  Lemma $3.3$ from Huang-Yin's paper ~\cite{HY2}):
\bl Let $N_{s}:=ts+2$. For all $0\leq j\leq t-1$ and
$p\in\left[2t+j(s-2)+2,2t+(j+1)(s-2)+1\right]$, we make the
following estimate
\begin{equation}\left.
\begin{split}G_{\geq p}(z,w)=& 2(1-s)^{j+1}\Re\left\{\langle z,a \rangle \Delta(z)^{j+1}
w^{t-j-1}\right\}+2\Re\left\langle \overline{z}+
\Delta'(z)+\Theta_{s}^{2}(z,\overline{z}),\overline{F_{\geq
p-1}(z,w)}\right\rangle\\&+\varphi'_{\geq
p}(z,\overline{z})-\varphi_{\geq
p}(z,\overline{z})+\Theta_{N_{s}}^{p}(z,\overline{z}),
\end{split}\right.\label{30}
\end{equation}
where
$\rm{wt}\left\{\overline{\Theta_{N_{s}}^{2t+2}(z,\overline{z})}\right\}\geq
 N_{s}$ and $w$ satisfies (\ref{ec4}).
  \el
\begin{proof}
 $\quad$

$\bf{Step\hspace{0.1 cm}1.}$ When $s=3$ this step is obvious.
Assume that $s>3$. Let  $p_{0}=2t+j(s-2)+2$, where
$j\in\left[0,t-1\right]$.
 We make induction on $p\in\left[2t+j(s-2)+2,2t+(j+1)(s-2)+1\right]$. For $j=0$
(therefore $p=2t+2$) the lemma is satisfied (see equation
(\ref{27})). Let $p\geq p_{0}$ such that $p+1\leq
2t+(j+1)(s-2)+1$. Collecting the terms of bidegree $(m,n)$ in
$(z,\overline{z})$ from (\ref{30}) with $m+n=p$, we obtain
\begin{equation}\hspace{0.1 cm}G_{nor}^{(p)}\left(z,\langle z,z \rangle\right)=
2\Re\left\langle z,F_{nor}^{(p-1)}\left(z,\langle
z,z\rangle\right)\right\rangle+
\varphi'_{p}(z,\overline{z})-\varphi_{p}(z,\overline{z})+\mathbb{P}_{N_{s}}^{p}(z,\overline{z}).\label{A}\end{equation}
By applying Extended Moser Lemma we  find a solution
$\left(F_{nor}^{(p-1)}(z,w),G_{nor}^{(p)}(z,w)\right)$ for
(\ref{A}). Assume that $p$ is even. In this case we  find
$F_{nor}^{(p-1)}(z,w)$ recalling the cases $1$ and $3$ of the
Extended Moser Lemma's proof. By using the cases $2$ and $4$ of the
Extended Moser Lemma's proof we find $G_{nor}^{(p)}(z,w)$. Since
$\mbox{wt}\left\{\overline{\mathbb{P}_{N_{s}}^{2t+2}(z,\overline{z})}\right\}\geq
 N_{s}$ we obtain $\mbox{wt}\left\{\left\langle
F_{nor}^{(p-1)}\left(z,\langle
z,z\rangle\right),z\right\rangle\right\},\mbox{wt}\left\{\left\langle
F_{nor}^{(p-1)}\left(z,\langle
z,z\rangle\right),z\right\rangle\right\}\geq N_{s}$.  Also
$\mbox{wt}\left\{G_{nor}^{(p)}\left(z,\langle
z,z\rangle\right)\right\},\hspace{0.1
cm}\mbox{wt}\left\{\overline{G_{nor}^{(p)}\left(z,\langle
z,z\rangle\right)}\right\}\geq N_{s}$. We can bring similarly
arguments as well when $p$ is even. We obtain the following
estimates
\begin{equation}\left.\begin{split}&\quad\quad\quad\quad\quad\quad\hspace{0.1 cm}\mbox{wt}\left\{F_{nor}^{(p-1)}(z,w)\right\}\geq N_{s}-s+1,\hspace{0.1
cm}\mbox{wt}\left\{\overline{F_{nor}^{(p-1)}\left(z,w\right)}\right\}\geq
N_{s}-1,\\&\mbox{wt}\left\{\overline{F_{nor}^{(p-1)}(z,w)-F_{nor}^{(p-1)}\left(z,\langle
z,z\rangle\right)}\right\}\geq N_{s}-1,
\quad\mbox{wt}\left\{F_{nor}^{(p-1)}(z,w)-F_{nor}^{(p-1)}\left(z,\langle
z,z\rangle\right)\right\}\geq N_{s}-s+1,
\\&\quad\quad\quad\quad\quad\quad\mbox{wt}\left\{G_{nor}^{(p)}(z,w)\right\}\geq
N_{s},\hspace{0.1
cm}\mbox{wt}\left\{G_{nor}^{(p)}(z,w)-G_{nor}^{(p)}\left(z,\langle
z,z\rangle\right)\right\}\geq N_{s},
\end{split}\right.\label{31345}\end{equation}
where  $w$ satisfies (\ref{ec4}). As a consequence of (\ref{A}) we
obtain
\begin{equation}\left.\begin{split}&G_{nor}^{(p)}(z,w)-G_{nor}^{(p)}\left(z,\langle
z,z\rangle\right)=\Theta_{N_{s}}^{p+1}(z,\overline{z})',\quad
2\Re\left\langle
z,F_{nor}^{(p-1)}(z,w)-F_{nor}^{(p-1)}\left(z,\langle
z,z\rangle\right)
\right\rangle=\Theta_{N_{s}}^{p+1}(z,\overline{z})'
,\\&\quad\quad\quad\left\langle\Delta'(z)+\Theta_{s}^{2}(z,\overline{z}),\overline{F_{nor}^{(p-1)}(z,w)}\right\rangle+\left\langle
\overline{F_{nor}^{(p-1)}(z,w)},\Delta'(z)+\Theta_{s}^{2}(z,\overline{z})\right\rangle=\Theta_{N_{s}}^{p+1}(z,\overline{z})',
\end{split}\right.\label{400}\end{equation}
and each of the preceding formal power series
$\Theta_{N_{s}}^{p+1}(z,\overline{z})'$ has the property
$\mbox{wt}\left\{\overline{\Theta_{N_{s}}^{p+1}(z,\overline{z})'}\right\}\geq
 N_{s}$. Substituting  $F_{\geq p-1}(z,w)=F_{nor}^{(p-1)}(z,w)+F_{\geq
p}(z,w)$ and $G_{\geq p}(z,w)=G_{nor}^{(p)}(z,w)+G_{\geq
p+1}(z,w)$ into (\ref{30}), we obtain
\begin{equation}\left.
\begin{split} G_{nor}^{(p)}(z,w)+G_{\geq p+1}(z,w)=&2(1-s)^{j+1}\Re\left\{\langle z,a\rangle\Delta(z)^{j+1}w^{t-j-1}\right\}
\\&+2\Re\left\langle\overline{z}+
\Delta'(z)+\Theta_{s}^{2}(z,\overline{z}),\overline{F_{nor}^{(p-1)}(z,w)}+\overline{F_{\geq
p}(z,w)}\right\rangle+\varphi'_{p}(z,\overline{z})-\varphi_{p}(z,\overline{z})\\&+\mathbb{P}_{N_{s}}^{p}(z,\overline{z})
+\varphi'_{\geq p+1}(z,\overline{z})-\varphi_{\geq
p+1}(z,\overline{z})+\Theta_{N_{s}}^{p+1}(z,\overline{z}).\end{split}\right.\label{32}\end{equation}

 Collecting the pure terms of degree $p$ from (\ref{A}), it follows that
$\varphi_{p,0}(z)=\varphi'_{p,0}(z)+\dots$, where in ''$\dots$''
we have determined terms with the weight less than
$p<N_{s}:=ts+2$. Therefore $\varphi_{p,0}(z)=\varphi'_{p,0}(z)$.
We will obtain that $\varphi_{k,0}(z)=\varphi'_{k,0}(z)$, for all
$k=3,\dots,T$. By making a simplification in (\ref{32}) by using
(\ref{A}),  it follows that
\begin{equation}\left.\begin{split} G_{\geq
p+1}(z,w)=&2(1-s)^{j+1}\Re\left\{\langle z,a
\rangle\Delta(z)^{j+1}w^{t-j-1}\right\}+2\Re\left\langle\overline{z}+
\Delta'(z)+\Theta_{s}^{2}(z,\overline{z}),\overline{F_{\geq
p}(z,w)}\right\rangle\\&+\varphi'_{\geq
p+1}(z,\overline{z})-\varphi_{\geq
p+1}(z,\overline{z})+J(z,\overline{z})+\Theta_{N_{s}}^{p+1}(z,\overline{z}),\end{split}\right.\label{33}
\end{equation}
where
$\mbox{wt}\left\{\overline{\Theta_{N_{s}}^{p+1}(z,\overline{z})}\right\}\geq
 N_{s}$ and
 \begin{equation}\left.\begin{split}
J(z,\overline{z})=&2\Re \left\langle
z,F_{nor}^{(p-1)}(z,w)-F_{nor}^{(p-1)}\left(z,\langle
z,z\rangle\right)\right\rangle+2\Re \left\langle
\Delta'(z)+\Theta_{s}^{2}(z,\overline{z}),
\overline{F_{nor}^{(p-1)}(z,w)}\right\rangle\\&+G_{nor}^{(p)}\left(z,\langle
z,z\rangle\right)-G_{nor}^{(p)}(z,w).
\end{split}\right.\label{VB}\end{equation}
By using (\ref{31345}) and (\ref{400}) it follows that
$J(z,\overline{z})=\Theta_{N_{s}}^{p+1}(z,\overline{z})$, where
$\mbox{wt}\left\{\overline{\Theta_{N_{s}}^{p+1}(z,\overline{z})}\right\}\geq
 N_{s}$.

 $\bf{Step\hspace{0.1 cm}2.}$ Assume that
we have proved Lemma $3.7$ for
$p\in\left[2t+j(s-2)+2,2t+(j+1)(s-2)+1\right]$ for $j\in[0,t-1]$.
We prove now  Lemma $3.7$ for
$p\in\left[2t+(j+1)(s-2)+2,2t+(j+2)(s-2)+1\right]$. Collecting the
terms of bidegree $(m,n)$ in $(z,\overline{z})$ from (\ref{30})
with $m+n=\Lambda+1:=2t+(j+1)(s-2)+1$,  we obtain
\begin{equation} \left.\begin{split} G_{nor}^{\left(\Lambda+1\right)}(z,\langle
z,z\rangle)=&2(1-s)^{j+1}\Re\left\{\langle
 z,a\rangle\Delta(z)^{j+1}\langle z,z\rangle^{t-j-1}\right\}+ 2\Re\left\langle
z,F_{nor}^{(\Lambda)}(z,\langle
z,z\rangle)\right\rangle\\&+\varphi'_{\Lambda+1}(z,\overline{z})-\varphi_{\Lambda+1}(z,\overline{z})
+\mathbb{P}_{N_{s}}^{\Lambda+1}(z,\overline{z}).
\end{split}\right.\label{34} \end{equation}

Here
$\mbox{wt}\left\{\overline{\left(\Theta_{1}\right)_{N_{s}}^{\Lambda+1}(z,\overline{z})}\right\}\geq
 N_{s}$. We define the following map
\begin{equation}F_{nor}^{(\Lambda)}(z,w)=F_{1}^{(\Lambda)}(z,w)+F_{2}^{(\Lambda)}(z,w),\quad F_{1}^{(\Lambda)}(z,w)=-(1-s)^{j+1}\langle z,a\rangle
\Delta(z)^{j+1}w^{t-j-2}\left(z_{1},\dots,z_{N}\right).\label{35}\end{equation}

Substituting (\ref{35})  into (\ref{34}), we obtain
\begin{equation}
G_{nor}^{\Lambda+1}\left(z,\langle
z,z\rangle\right)=2\Re\left\langle z,F_{2}^{(\Lambda)}(z,\langle
z,z\rangle)\right\rangle+
\varphi'_{\Lambda+1}(z,\overline{z})-\varphi_{\Lambda+1}(z,\overline{z})+\mathbb{P}_{N_{s}}^{\Lambda+1}(z,\overline{z}).\label{36}
\end{equation}
By applying Extended Moser Lemma we  find a solution
$\left(G_{nor}^{(\Lambda+1)}(z,w),F_{2}^{(\Lambda)}(z,w)\right)$
for
 (\ref{36}). By using the same  arguments as in the Step $1$ we obtain the following
estimates
\begin{equation}\left.\begin{split}&
\mbox{wt}\left\{G_{nor}^{(\Lambda+1)}(z,w)-G_{nor}^{(\Lambda+1)}\left(z,\langle
z,z\rangle\right)\right\},\mbox{wt}\left\{G_{nor}^{(\Lambda+1)}(z,w)\right\},\hspace{0.1
cm}\mbox{wt}\left\{ G_{nor}^{(\Lambda+1)}\left(z,\langle
z,z\rangle\right) \right\} \geq N_{s}
,\\&\mbox{wt}\left\{F_{2}^{(\Lambda)}(z,w)-F_{2}^{(\Lambda)}\left(z,\langle
z,z\rangle\right)\right\},\hspace{0.2 cm}
\mbox{wt}\left\{F_{2}^{(\Lambda)}(z,w)\right\},
\mbox{wt}\left\{F_{2,k}^{(\Lambda)}\left(z,\langle
z,z\rangle\right)\right\}\geq N_{s}-s+1,\\&\quad
\mbox{wt}\left\{\overline{
F_{2}^{(\Lambda)}(z,w)-F_{2}^{(\Lambda)}\left(z,\langle z,z
\rangle\right)}\right\},  \mbox{wt}\left\{
\overline{F_{2}^{(\Lambda)}(z,w)}\right\},\hspace{0.1
cm}\mbox{wt}\left\{\overline{F_{2}^{(\Lambda)}\left(z,\langle
z,z\rangle\right)}\right\}\geq N_{s}-1,
\end{split}\right.\label{31}\end{equation} where  $w$ satisfies
(\ref{ec4}). As a consequence of (\ref{31}) we obtain
\begin{equation}\left.\begin{split}&\quad\quad
\quad\quad\left\langle\Delta'(z)+\Theta_{s}^{2}(z,\overline{z}),\overline{F_{2}^{(\Lambda)}(z,w)}\right\rangle+\left
\langle
\overline{F_{2}^{(\Lambda)}(z,w)},\Delta'(z)+\Theta_{s}^{2}(z,\overline{z})\right\rangle=\Theta_{N_{s}}^{\Lambda+2}(z,\overline{z})',
\\&
G_{nor}^{(\Lambda+1)}(z,w)-G_{nor}^{(\Lambda+1)}\left(z,\langle
z,z\rangle\right)=\Theta_{N_{s}}^{\Lambda+2}(z,\overline{z})',
\quad2\Re\left\langle
F_{2}^{(\Lambda)}(z,w)-F_{2}^{(\Lambda)}\left(z,\langle
z,z\rangle\right),z\right\rangle=\Theta_{N_{s}}^{\Lambda+2}(z,\overline{z})',
\end{split}\right.\label{37}
\end{equation}
where $w$ satisfies (\ref{ec4}) and each of the preceding formal
power series has the property
$\mbox{wt}\left\{\overline{\Theta_{N_{s}}^{\Lambda+2}(z,\overline{z})}\right\}\geq
 N_{s}$. Substituting
$F_{\geq\Lambda}(z,w)=F_{nor}^{(\Lambda)}(z,w)+F_{\geq\Lambda+1}(z,w)$
and
$G_{\geq\Lambda+1}(z,w)=G_{nor}^{(\Lambda+1)}(z,w)+G_{\geq\Lambda+2}(z,w)$
in (\ref{30}), we obtain
\begin{equation}\left. \begin{split} G_{nor}^{(\Lambda+1)}(z,w)+G_{\geq\Lambda+2}(z,w)&=2(1-s)^{j+1}\Re\left\{\langle z,a
\rangle\Delta(z)^{j+1}w^{t-j-1}\right\}\\&\quad+2\Re\left\langle
\overline{z}+ \Delta'(z)
+\Theta_{s}^{2}(z,\overline{z}),\overline{F_{nor}^{(\Lambda)}(z,w)}+\overline{F_{\geq
\Lambda+1}(z,w)}\right\rangle
+\varphi_{\Lambda+1}(z,\overline{z})-\varphi'_{\Lambda+1}(z,\overline{z})
\\&\quad+\varphi'_{\geq\Lambda+2}(z,\overline{z})-\varphi_{\geq\Lambda+2}(z,\overline{z})+\mathbb{P}_{N_{s}}^{\Lambda+1}(z,\overline{z})
+\Theta_{N_{s}}^{\Lambda+2}(z,\overline{z}).\end{split}\right.\label{38}\end{equation}
By making a simplification in (\ref{38}) with (\ref{34}), and then
by using (\ref{35}), we obtain
\begin{equation}\left.\begin{split}G_{\geq\Lambda+2}(z,w)=
2\Re\left\langle \overline{z}+
\Delta'(z)+\Theta_{s}^{2}(z,\overline{z}),\overline{F_{\geq\Lambda+1}(z,w)}\right\rangle
+\varphi'_{\geq\Lambda+2}(z,\overline{z})-\varphi_{\geq\Lambda+2}(z,\overline{z})
+\Theta_{N_{s}}^{\Lambda+2}(z,\overline{z})+J(z,\overline{z}),\end{split}\right.\label{39}\end{equation}
where
\begin{equation}\left.\begin{split}J(z,\overline{z})&=2\Re\left\langle z,F_{nor}^{(\Lambda)}(z,w)-F_{nor}^{(\Lambda)}\left(z,\langle
z,z\rangle\right)\right\rangle+2\Re\left\langle\Delta'(z)+\Theta_{s}^{2}(z,\overline{z}),\overline{F_{nor}^{(\Lambda)}(z,w)}\right\rangle\\&\quad
+2(1-s)^{j+1}\Re\left\{\langle
z,a\rangle\Delta(z)^{j+1}w^{t-j-1}-\langle z,a\rangle
\Delta(z)^{j+1}\langle
z,z\rangle^{t-j-1}\right\}\\&\quad+G_{nor}^{(\Lambda)}\left(z,\langle
z,z\rangle\right)-G_{nor}^{(\Lambda)}(z,w),\\&=2\Re\left\langle
z,F_{1}^{(\Lambda)}(z,w) -F_{1}^{(\Lambda)}\left(z,\langle
z,z\rangle\right)+F_{2}^{(\Lambda)}(z,w)-F_{2}^{(\Lambda)}\left(z,\langle
z,z\rangle\right)\right\rangle\\&\quad+2\Re\left\langle\Delta'(z)+\Theta_{s}^{2}(z,\overline{z})
,\overline{F_{1}^{(\Lambda)}(z,w)}+\overline{F_{2}^{(\Lambda)}(z,w)}\right\rangle+G_{nor}^{(\Lambda)}\left(z,\langle
z,z\rangle\right)-G_{nor}^{(\Lambda)}(z,w)\\&\quad+2(1-s)^{j+1}\Re\left\{\langle
z,a\rangle\Delta(z)^{j+1}\left(w^{t-j-1}-\langle
z,z\rangle^{t-j-1}\right)\right\}.\end{split}\right.\label{40}\end{equation}
By using (\ref{31}) and (\ref{37}) it follows that
\begin{equation}\left.\begin{split}J(z,\overline{z})&
= 2\Re\left\langle
z,F_{1}^{(\Lambda)}(z,w)-F_{1}^{(\Lambda)}\left(z,\langle
z,z\rangle\right)\right\rangle+2\Re\left\langle
\Delta'(z)+\Theta_{s}^{2}(z,\overline{z})
,\overline{F_{1}^{(\Lambda)}(z,w)}\right\rangle\\&\quad+2(1-s)^{j+1}\Re\left\{\langle
z,a\rangle \Delta(z)^{j+1}\left(w^{t-j-1}-\langle
z,z\rangle^{t-j-1}\right)\right\}+\Theta_{N_{s}}^{\Lambda+2}(z,\overline{z}),\end{split}\right.\label{100}\end{equation}
where
$\mbox{wt}\left\{\overline{\Theta_{N_{s}}^{\Lambda+2}(z,\overline{z})}\right\}\geq
 N_{s}$. We observe that
\begin{equation}\Re \left\langle
z,F_{1}^{(\Lambda)}\left(z,\langle
z,z\rangle\right)\right\rangle=-(1-s)^{j+1}\Re\left\{\langle
z,a\rangle\langle
z,z\rangle^{t-j-1}\Delta(z)^{j+1}\right\}.\label{900}
\end{equation}
 Since
$\mbox{wt}\left\{F_{1}^{\left(\Lambda\right)}(z,w)\right\}\geq
N_{s}-s$ and
 $\mbox{wt}\left\{\overline{F_{1}^{(\Lambda)}(z,w)}\right\}\geq N_{s}$, it follows that
\begin{equation}\Re\left\langle \Theta_{s}^{2}(z,\overline{z}),
\overline{F^{\left(\Lambda\right)}_{1}(z,w)}\right\rangle=\Theta_{N_{s}}^{\Lambda+2}(z,\overline{z}),
\label{44}\end{equation} where
$\mbox{wt}\left\{\overline{\Theta_{N_{s}}^{\Lambda+2}(z,\overline{z})}\right\}\geq
 N_{s}$. By using (\ref{900}) and (\ref{44}),
we can rewrite (\ref{100}) as follows
\begin{equation}\left.\begin{split}J(z,\overline{z})
=&2\Re\left\langle
z,F_{1}^{\left(\Lambda\right)}(z,w)\right\rangle+2\Re\left\langle
\Delta'(z),\overline{F_{1}^{\left(\Lambda\right)}(z,w)}\right\rangle+2(1-s)^{j+1}\Re\left\{\langle
z,a\rangle\Delta(z)^{j+1}
w^{t-j-1}\right\}+\Theta_{N_{s}}^{\Lambda+2}(z,\overline{z}),\end{split}\right.\label{101}\end{equation}
where
$\mbox{wt}\left\{\overline{\Theta_{N_{s}}^{\Lambda+2}(z,\overline{z})}\right\}\geq
 N_{s}$.  Substituting the
formula of $F_{1}^{\left(\Lambda\right)}(z,w)$ in (\ref{101}), we
obtain
\begin{equation}\left.\begin{split}J(z,\overline{z})
&=-2(1-s)^{j+1}\Re\left\{\langle
z,a\rangle\Delta(z)^{j+1}w^{t-j-2}\left(\langle
z,z\rangle+\left\langle
\left(z_{1},\dots,z_{N}\right),\left(\overline{\Delta_{1}(z)},\dots,\overline{\Delta_{N}(z)}\right)\right\rangle\right)\right\}
\\&\quad+2(1-s)^{j+1}\Re\left\{\langle
z,a\rangle\Delta(z)^{j+1}
w^{t-j-1}\right\}+\Theta_{N_{s}}^{\Lambda+2}(z,\overline{z}),\\&=-2(1-s)^{j+1}\Re\left\{\langle
z,a\rangle\Delta(z)^{j+1}w^{t-j-2}\left(\langle
z,z\rangle+s\Delta(z)-w\right)\right\}+\Theta_{N_{s}}^{\Lambda+2}(z,\overline{z}),\\&=2(1-s)^{j+2}\Re\left\{\langle
z,a\rangle\Delta(z)^{j+2}w^{t-j-2}\right\}+\Theta_{N_{s}}^{\Lambda+2}(z,\overline{z}),
\end{split}\right.\label{102}\end{equation}
where $w$ satisfies (\ref{ec4}) and
$\mbox{wt}\left\{\overline{\Theta_{N_{s}}^{\Lambda+2}(z,\overline{z})}\right\}\geq
 N_{s}$.

 The proof of our Lemma follows by using (\ref{102}) and (\ref{39}).
 \end{proof}

Collecting the  terms of bidegree $(m,n)$ in $(z,\overline{z})$
with $m+n=ts+1$ and $t=j-1$ from (\ref{30}), we obtain
\begin{equation}\left.\begin{split} G_{nor}^{(ts+1)}\left(z,\langle z,z\rangle\right)=&2(1-s)^{t}\Re\left\{\langle z,a \rangle
\Delta(z)^{t}\right\} +2\Re\left\langle
z,F_{nor}^{(ts)}\left(z,\langle z,z\rangle\right)\right\rangle\\&
+\varphi'_{ts+1,0}(z,\overline{z})-\varphi_{ts+1,0}(z,\overline{z})+\left(\Theta_{1}\right)_{N_{s}}^{ts+1}(z,\overline{z}).\end{split}\right.\label{46}\end{equation}
By applying Extended Moser Lemma we find a solution
$\left(G_{nor}^{(ts+1)}(z,w),F_{nor}^{(ts)}(z,w)\right)$ for
(\ref{46}). Collecting the pure terms from (\ref{46}) of degree
$ts+1$, it follows that
\begin{equation}\varphi'_{ts+1,0}(z)-\varphi_{ts+1,0}(z)=(1-s)^{t}\langle
z,a\rangle\Delta(z)^{t}.\label{caz1}\end{equation} The parameter
$a$ will help us to put the desired normalization condition (see
(\ref{cnfnp})). By applying Lemma $2.4$ for $\varphi'_{ts+1,0}(z)$
and $\varphi_{ts+1,0}(z)$, it follows that
\begin{equation}\varphi_{ts+1,0}(z)=(1-s)^{t}Q(z)\Delta(z)^{t}+R(z),\quad
\varphi'_{ts+1,0}(z)=Q'(z)\Delta(z)^{t}+R'(z),\label{47}\end{equation}
where
$\left(\Delta^{t}\right)^{\star}\left(R(z)\right)=\left(\Delta^{t}\right)^{\star}\left(R'(z)\right)=0$.
We impose the normalization condition
 $\left(\Delta^{t}\right)^{\star}\left(\varphi'_{ts+1,0}(z)\right)=0$. This is equivalent finding $a$ such that
 $Q'(z)=0$. Here $Q(z)$ is a determined holomorphic polynomial. We find $a$ by solving the equation
$Q'(z)=(1-s)^{t}\langle z,a\rangle-Q(z)=0.$

\vspace{0.1 cm}

 By composing the map that sends $M$ into
(\ref{ec4}) with the map (\ref{t1}) we obtain our formal
transformation that sends $M$ into $M'$ up  to degree $ts+1$.

\section{Proof of Theorem $1.3$-Case $T+1=(t+1)s$, $t\geq1$}

 In this case we are looking for a biholomorphic transformation of the following
type
\begin{equation}\left.\begin{split}&\quad\quad\quad\quad\quad\left(z',w'\right)=\left(z+F(z,w),w+G(z,w)\right)\\ &F(z,w)=\displaystyle\sum_{l=0}^{T-2t-1}F_{nor}^{(2t+l+1)}(z,w),
\hspace{0.1cm}G(z,w)=\displaystyle\sum_{\tau=0}^{T-2t}G_{nor}^{(2t+2+\tau)}(z,w)\end{split}\right.\label{48},\end{equation}
that maps $M$ into $M'$ up to the degree $T+1=(t+1)s$. In order to
make the  mapping (\ref{48})  uniquely determined we assume that
$F_{nor}^{(2t+l+1)}(z,w)$ is normalized as in Extended Moser
Lemma, for all $l=1,\dots T-2t-1$. Replacing  (\ref{48}) in
(\ref{em2}), and after a simplification with (\ref{ec4}),
 we obtain
\begin{equation}\left.\begin{split}\displaystyle &\sum_{\tau=0}^{T-2t-1}G_{nor}^{(2t+2+\tau)}\left(z,\langle z,z\rangle+\varphi_{\geq 3}(z,\overline{z})\right)
=2\Re\left\langle
\displaystyle\sum_{l=0}^{T-2t-1}F_{nor}^{(2t+l+1)}\left(z,\langle
z,z\rangle+\varphi_{\geq 3}(z,\overline{z})\right), z
\right\rangle
\\&+\left\|\displaystyle\sum_{l=0}^{T-2t-1}F_{nor}^{(2t+l+1)}\left(z,\langle
z,z\rangle+\varphi_{\geq 3}(z,\overline{z})\right)\right\|^{2}\\&+
\varphi'_{\geq
3}\left(z+\displaystyle\sum_{l=-1}^{T-2t}F_{nor}^{(2t+l+2)}\left(z,\langle
z,z\rangle+\varphi_{\geq
3}(z,\overline{z})\right),\overline{z+\displaystyle\sum_{l=-1}^{T-2t}F_{nor}^{(2t+l+2)}\left(z,\langle
z,z\rangle+\varphi_{\geq
3}(z,\overline{z})\right)}\right)-\varphi_{\geq
3}(z,\overline{z}).\end{split}\right.\label{49}\end{equation}

Collecting the terms with the same bidegree in $(z,\overline{z})$
from (\ref{49})  we will find  $F(z,w)$ and $G(z,w)$ by applying
Extended Moser Lemma. Since $F(z,w)$ and $G(z,w)$ don't have
components of normal weight less than $2t+2$, collecting in
(\ref{49}) the terms of bidegree $(m,n)$ in $(z,\overline{z})$
with $m+n<2t+2$, we obtain
$\varphi'_{m,n}(z,\overline{z})=\varphi_{m,n}(z,\overline{z})$.

Collecting the  terms of bidegree $(m,n)$  in $(z,\overline{z})$
with $m+n=2t+2$ from (\ref{49}), we prove the following lemma:
 \bl $G_{nor}^{(2t+2)}(z,w)=\left(a+\overline{a}\right)w^{t+1},\hspace{0.1
cm} F_{nor}^{(2t+1)}(z,w)=w^{t}\left(\begin{array}{ccc}
a_{11}&\dots&
a_{1N}\\ \vdots & \ddots & \vdots \\
a_{N1}& \dots & a_{NN}\end{array}\right)\left(\begin{array}{l}
z_{1}\\ \vdots \\
z_{N}\end{array}\right)$, where $Na$ is the trace of the matrix
$\left(a_{ij}\right)_{1\leq i,j \leq N}$. \el
\begin{proof}

Collecting the pure terms of degree $2t+2$  from (\ref{49}), we
obtain that $\varphi_{2t+2}(z)=\varphi'_{2t+2}(z)$. Collecting the
 terms of bidegree $(m,n)$ in $(z,\overline{z})$ with
$m+n=2t+2$ and $0<m<n-1$ from (\ref{49}), we obtain
\begin{equation}\varphi_{m,n}'(z,\overline{z})=-\left\langle
z,F_{n-m+1,m-1}(z)\right\rangle \langle
z,z\rangle^{m-1}+\varphi_{m,n}(z,\overline{z}).\end{equation}
Since $\varphi_{m,n}(z,\overline{z}),\hspace{0.1
cm}\varphi'_{m,n}(z,\overline{z})$ satisfy (\ref{CN1}), by the
uniqueness of trace decomposition, we obtain $F_{n-m+1,m-1}(z)=0$.
Collecting the  terms of bidegree $(m,n)$ in $(z,\overline{z})$
with $m+n=2t+2$ and $m>n+1$ from (\ref{49}), we obtain
\begin{equation}\varphi_{m,n}'(z,\overline{z})=G_{m-n}(z)\langle
z,z\rangle^{n}-\left\langle F_{m-n+1,n-1}(z),z\right\rangle
\langle
z,z\rangle^{n-1}+\varphi_{m,n}(z,\overline{z}).\end{equation}
Since $F_{n-m+1,m-1}(z)=0$ it follows that $G_{m-n}(z)=0$.

Collecting the terms of bidegree $(t+1,t+1)$ in $(z,\overline{z})$
from (\ref{49}), we obtain
\begin{equation}\varphi_{t+1,t+1}'(z,\overline{z})=\left(G_{0,t+1}(z)\langle
z,z\rangle-\left\langle F_{1,t}(z),z\right \rangle-\left\langle
z,F_{1,t}(z)\right\rangle\right)\langle
z,z\rangle^{t}+\varphi_{t+1,t+1}(z,\overline{z}).\label{908}
\end{equation} Then (\ref{908}) can not provide us $F_{1,t}(z)$.
Therefore $F_{1,t}(z)$ is undetermined. We obtain
\begin{equation}F_{nor}^{(2t+1)}(z,w)=w^{t}\left(\begin{array}{ccc}
a_{11}&\dots&a_{1N}\\ \vdots&\ddots&\vdots\\a_{N1}&\dots&a_{NN}
\end{array}\right)
\left(\begin{array}{l} z_{1}\\\vdots\\
z_{N}\end{array}\right),\quad a_{ij}\in\mathbb{C},\hspace{0.1 cm}
1\leq i,j \leq N.\end{equation} We  write
$a_{11}=a+b_{11},\dots,a_{NN}=a+b_{NN}$ and we  use the notations
$b_{k,j}=a_{k,j}$, for all $k\neq j$. Then the matrix
$\left(b_{k,j}\right)_{1\leq k,j\leq N}$ represents the traceless
part of the matrix $\left(a_{k,j}\right)_{1\leq k,j\leq N}$. By
applying Lemma $2.1$ to the polynomial $\left\langle
F_{1,t}(z),z\right\rangle$, we obtain $\left\langle
F_{1,t}(z),z\right\rangle=a\langle z,z\rangle+P(z,\overline{z})$
with $\mbox{tr}\left(P(z,\overline{z})\right)=0$,  where
$P(z,\overline{z})=\displaystyle\sum_{i,j=1}^{N}b_{i,j}z_{i}\overline{z}_{j}$.
By using the preceding decomposition we obtain
   \begin{equation}\varphi_{t+1,t+1}'(z,\overline{z})=\left(G_{0,t+1}(z)-a-\overline{a}\right)\langle z,z \rangle^{t+1}
     +\varphi_{t+1,t+1}(z,\overline{z})-2\Re\left(P(z,\overline{z})\langle
     z,z\rangle^{t}\right).\label{GR}\end{equation}
Since $\mbox{tr}\left(P(z,\overline{z})\right)=0$ it follows that
$\mbox{tr}^{t+1}\left(\Re \left(P(z,\overline{z})\langle
z,z\rangle^{t}\right)\right)=0$ (see Lemma $6.6$ from
~\cite{DZ1}).
 \end{proof}

We can write $F(z,w)=F_{nor}^{(2t+2)}(z,w)+F_{\geq2t+3}(z,w)$ and
$G(z,w)=G_{\geq2t+2}(z,w)$ (see (\ref{2.11})). We have $F_{\geq
2t+2}(z,w)=\displaystyle\sum_{k+2l\geq 2t+2}F_{k,l}(z)w^{l}$,
where $F_{k,l}(z)$ is a homogeneous polynomial of degree $k$.
Therefore $\mbox{wt}\left\{F_{\geq
2t+2}(z,w)\right\}\geq\displaystyle\min_{k+2l\geq2t+2}\{k+ls\}\geq\displaystyle\min_{k+2l\geq2t+2}\{k+2l\}\geq
2t+2.$ Next, we show that
$\mbox{wt}\left\{\overline{F_{\geq2t+2}(z,w)}\right\}\geq ts+s-1$.
 Since  $\mbox{wt}\left\{\overline{F_{\geq2t+2}(z,w)}\right\}\geq\displaystyle\min_{k+2l\geq2t+2}\{k(s-1)+ls\}$,
it is enough to prove that $k(s-1)+ls\geq ts+s-1$ for $k+2l\geq
2t+2$. Since we can write the latter inequality as
$(k-1)(s-1)+ls\geq ts$ for $(k-1)+2l\geq 2t+1$, it is enough to
prove that $k(s-1)+ls\geq ts$ for $k+2l\geq 2t+1>2t$. Continuing
the calculations like in the previous case we obtain the desired
result.

\bl For $w$ satisfying (\ref{ec4}), we make the following
immediate estimates
\begin{equation}\left.\begin{split}&\rm{wt}\left\{F_{nor}^{(2t+1)}(z,w)\right\}\geq
ts+1,\quad
\rm{wt}\left\{\overline{F_{nor}^{(2t+1)}(z,w)}\right\}\geq
ts+s-1,\quad\rm{wt}\left\{\left\|F_{nor}^{(2t+1)}(z,w)\right\|^{2}\right\}\geq
ts+s+1,\\&\quad\rm{wt}\left\{F_{\geq 2t+2}(z,w)\right\}\geq
2t+2,\quad \left\{\overline{F_{\geq 2t+2}(z,w)} \right\}\geq
ts+s-1,
\quad\rm{wt}\left\{\left\|F_{\geq2t+2}(z,w)\right\|^{2}\right\}\geq
ts+s+1,\\&\quad\quad\quad\hspace{0.1 cm}\rm{wt}\left\{\left\langle
F_{nor}^{(2t+1)}(z,w),F_{\geq
2t+2}(z,w)\right\rangle\right\},\hspace{0.1
cm}\mbox{wt}\left\{\left\langle F_{\geq 2t+2}(z,w),
F_{nor}^{(2t+1)}(z,w)\right\rangle\right\}\geq ts+s+1.
\end{split}\right.\label{50}\end{equation}

\el As a consequence of the estimates (\ref{50}) we obtain
\begin{equation}\left\|F(z,w)\right\|^{2}=\left\|F_{nor}^{(2t+1)}(z,w)\right\|^{2}+2\Re\left\langle
F_{nor}^{(2t+1)}(z,w),F_{\geq 2t+2}(z,w)\right\rangle
+\left\|F_{\geq2t+2}(z,w)\right\|^{2}=\Theta_{ts+s+1}^{2t+3}(z,\overline{z})\label{500},\end{equation}
where
$\rm{wt}\left\{\overline{\Theta_{ts+s+1}^{2t+3}(z,\overline{z})}\right\}\geq
 ts+s+1$.

In order to apply Extended Moser Lemma in (\ref{49}) we have to
identify and weight and order evaluate the terms which are not
''good''. We prove the following lemmas:
 \bl For all $m,n\geq 1$ and $w$ satisfying (\ref{ec4}), we make the following
estimate
\begin{equation}\varphi'_{m,n}\left(z+F(z,w),\overline{z+F(z,w)}\right)
=\varphi'_{m,n}(z,\overline{z})+2\Re\left\langle
\Theta_{s}^{2}(z,\overline{z}),\overline{ F_{\geq
2t+2}(z,w)}\right\rangle+\Theta_{ts+s+1}^{2t+3}(z,\overline{z}),\label{901}\end{equation}
where
$\rm{wt}\left\{\overline{\Theta_{ts+s+1}^{2t+3}(z,\overline{z})}\right\}\geq
 ts+s+1$.
 \el
\begin{proof}

We have the  expansion
$\varphi'_{m,n}\left(z+F(z,w),\overline{z+F(z,w)}\right)
=\varphi'_{m,n}(z,\overline{z})+\dots$ (see the proof of  Lemma
$3.3$). In order to prove (\ref{901}), it is enough to study the
weight and the order of the following particular terms
$$A_{1}(z,w)=F_{1}(z,w)z^{I}\overline{z}^{J},\quad A_{2}(z,w)=
z^{I_{1}}\overline{z}^{J_{1}}\overline{F_{1}(z,w)} ,\quad
B_{1}(z,w)=F_{2}(z,w)z^{I}\overline{z}^{J},\quad
B_{2}(z,w)=z^{I_{1}}\overline{z}^{J_{1}}\overline{F_{2}(z,w)},$$
where $F_{1}(z,w)$ is the first component of
$F_{nor}^{(2t+1)}(z,w)$ and $F_{2}(z,w)$ is the first component of
$F_{\geq 2t+2}(z,w)$. Here we assume  that $\left|I\right|=m-1$,
$\left|J\right|=n$, $\left|I_{1}\right|=m$,
$\left|J_{1}\right|=n-1$.

By using (\ref{50}) we obtain
$\mbox{wt}\left\{A_{1}(z,w)\right\}\geq m-1+ts+1+n(s-1)\geq ts+s+1
\Longleftrightarrow m+ns-n\geq s+1\Longleftrightarrow m+s(n-1)\geq
n+1 $ and the latter inequality is true since $m+s(n-1)\geq
m+3(n-1)\geq n+1$. On the other hand
$\mbox{Ord}\left\{A_{1}(z,w)\right\}\geq m-1+2t+1+n\geq 2t+3.$

By using (\ref{50}) we obtain
$\mbox{wt}\left\{A_{2}(z,w)\right\}\geq m+(n-1)(s-1)+ts+s-1\geq
ts+s+1$ and the last inequality is equivalent with
$m+(n-1)(s-1)\geq 2$. The latter inequality can be proved with the
same calculations like in Lemma $3.3$ proof. On the other hand, we
observe that $\mbox{Ord}\left\{A_{1}(z,w)\right\}\geq
m+2t+1+n-1\geq 2t+3$.

In the same  way we obtain $\mbox{Ord}\left\{B_{1}(z,w)\right\}$,
$\mbox{Ord}\left\{B_{2}(z,w)\right\}\geq 2t+2$. By using
(\ref{50}), every term from ''$\dots$'' that depends on
$F_{2}(z,w)$ can be written as
$\Theta_{s}^{2}(z,\overline{z})F_{2}(z,w)$. This proves our claim.
\end{proof}
\bl For all $k>s$ and $w$ satisfying (\ref{ec4}), we make the
following estimate
\begin{equation}\varphi'_{k,0}\left(z+F(z,w)\right)=\varphi'_{k,0}(z)+2\Re\left\langle
\Theta_{s}^{2}(z,\overline{z}),\overline{ F_{\geq
2t+2}(z,w)}\right\rangle+\Theta_{ts+s+1}^{2t+3}(z,\overline{z}),\end{equation}
where
$\rm{wt}\left\{\overline{\Theta_{ts+s+1}^{2t+3}(z,\overline{z})}\right\}\geq
 ts+s+1$.
 \el
 \begin{proof}

We make the expansion
$\varphi'_{k,0}\left(z+F(z,w)\right)=\varphi'_{k,0}(z)+\dots$ . To
study the weight and the order  of terms which can appear in
''$\dots$'' it is enough to study the weight and order of the
following terms
$$A(z,w)=F_{1}(z,w)z^{I},\quad B(z,w)=F_{2}z^{I},$$ where $F_{1}(z,w)$ is the
first component of $F_{nor}^{(2t+1)}(z,w)$ and $F_{2}(z,w)$ is the
first component of $F_{\geq 2t+3}(z,w)$. Here we assume  that
$\left|I\right|=m-1\geq s$. From (\ref{50}) we obtain
 $\mbox{wt}\left\{A(z,w)\right\}\geq s+ts+1=ts+s+1$. On the other
hand, we have $\mbox{Ord}\left\{A(z,w)\right\}\geq 2t+s+1\geq
2t+3$. By using  (\ref{50}) each term from ''$\dots$'' that
depends on $F_{2}(z,w)$ can be written as
$\Theta_{s}^{2}(z,\overline{z})F_{2}(z,w)$. This proves our claim.
\end{proof}

\bl For $w$ satisfying   (\ref{ec4}) we have the following
estimate
\begin{equation}\left.\begin{split}2\Re\{\Delta\left(z+F(z,w)\right)\}
=&2\Re\left\{\Delta(z)+\displaystyle\sum_{k=1}^{N}\Delta_{k}(z)\left(a_{k1}z_{1}+\dots+a_{kN}z_{N}\right)w^{t}\right\}\\&+2\Re\left\langle
\Delta'(z)+\Theta_{s}^{2}(z,\overline{z}) ,\overline{F_{\geq
2t+2}(z,w)}\right\rangle+\Theta_{ts+s+1}^{2t+3}(z,\overline{z}),\end{split}\right.\label{24}\end{equation}
where
$\rm{wt}\left\{\overline{\Theta_{ts+s+1}^{2t+3}(z,\overline{z})}\right\}\geq
 ts+s+1$. \el
\begin{proof}

 For $w$ satisfying (\ref{ec4}), we have the expansion
\begin{equation}\left.\begin{split}2\Re\left\{\Delta\left(z+F(z,w)\right)\right\}=2\Re\left\{\Delta(z)+
\displaystyle\sum_{k=1}^{N}\Delta_{k}(z)F_{\geq
2t+1}^{k}(z,w)+L(z,\overline{z})\right\}+\Theta_{ts+s+1}^{2t+3}(z,\overline{z}),
\end{split}\right.\label{25}\end{equation}
where $F_{\geq 2t+1}(z,w)=\left(F^{1}_{\geq
2t+1}(z,w),\dots,F^{N}_{\geq 2t+1}(z,w)\right)$ and
$L(z,\overline{z})=\left\langle \Theta_{s}^{2}(z,\overline{z}),
\overline{F_{\geq2t+2}(z,w)}\right\rangle$. We compute
\begin{equation}\left.\begin{split} \displaystyle\sum_{k=1}^{N}2\Re\left\{\Delta_{k}(z)F_{\geq 2t+1}^{k}(z,w)\right\}&
=\displaystyle
\displaystyle\sum_{k=1}^{N}2\Re\left\{\Delta_{k}(z)\left(w^{t}\displaystyle\sum_{j=1}^{N}a_{kj}z_{j}+F_{\geq2t+2}^{k}(z,w)\right)\right\}
 \\&=2\Re\left\{
w^{t}\displaystyle\sum_{k=1}^{N}\Delta_{k}(z)\left(a_{k1}z_{1}+\dots+a_{kN}z_{N}\right)\right\}+2\Re\left\langle
\Delta'(z),\overline{F_{\geq 2t+2}(z,w)}\right\rangle.
\end{split}\right.\label{26}\end{equation}
 \end{proof}
 \bl For $w$ satisfying (\ref{ec4}), we have the following estimate
 \begin{equation}\left.\begin{split}
 G_{nor}^{(2t+2)}(z,w)-2\Re\left\langle
F_{nor}^{(2t+1)}(z,w),z \right\rangle
=&2(a+\overline{a})\Re\left\{\Delta(z)
w^{t}\right\}+2\Re\left\{P(z,\overline{z})w^{t}\right\}
+\Theta_{ts+s+1}^{2t+3}(z,\overline{z}),\end{split}\right.\end{equation}
where
$P(z,\overline{z})=\displaystyle\sum_{k,j=1}^{N}b_{k,j}z_{k}\overline{z}_{j}$
and
$\rm{wt}\left\{\overline{\Theta_{ts+s+1}^{2t+3}(z,\overline{z})}\right\}\geq
 ts+s+1$.
 \el
\begin{proof} For $w$ satisfying (\ref{ec4}), by Lemma $4.1$ it follows that
\begin{equation}\left.\begin{split}G_{nor}^{(2t+2)}(z,w)-2\Re\left\langle F_{nor}^{(2t+1)}(z,w),z \right\rangle  &=(a+\overline{a})w^{t+1}-2\Re
\Bigg\langle w^{t}\left(\begin{array}{ccc} b_{11}+a&\dots& a_{1N}\\
\vdots & \ddots & \vdots \\a_{N1} & \dots
&b_{NN}+a\end{array}\right)\left(\begin{array}{l} z_{1}\\ \vdots \\
z_{N}\end{array}\right),z \Bigg \rangle ,\\&=2\Re
\left\{aw^{t+1}\right\}-2\Re\left\{ aw^{t}\langle z,z\rangle
+P(z,\overline{z})w^{t}\right\}+\overline{a}\left(w^{t+1}-\overline{w}^{t+1}\right),
\\&=2\Re\left\{aw^{t}\left(w-\langle
z,z\rangle\right)\right\}-2\Re\left\{P(z,\overline{z})w^{t}\right\}+\Theta_{ts+s+1}^{2t+3}(z,\overline{z}),
\\&=2\Re\left\{aw^{t}\left(\Delta(z)+\overline{\Delta(z)}\right)\right\}-2\Re\left\{P(z,\overline{z})w^{t}\right\}+\Theta_{ts+s+1}^{2t+3}(z,\overline{z}),\\&
=2(a+\overline{a})\Re\left\{\Delta(z)w^{t}\right\}-2\Re\left\{P(z,\overline{z})w^{t}\right\}+\Theta_{ts+s+1}^{2t+3}(z,\overline{z}),
\end{split}\right.\label{52}\end{equation}
where
$\rm{wt}\left\{\overline{\Theta_{ts+s+1}^{2t+3}(z,\overline{z})}\right\}\geq
 ts+s+1$.
\end{proof}
Substituting  $F(z,w)=F^{(2t+1)}_{nor}(z,w)+F_{\geq
 2t+2}(z,w)$ and $G(z,w)=G^{(2t+2)}_{nor}(z,w)+G_{\geq 2t+3}(z,w)$ (see (\ref{2.11})) into (\ref{49}) and by Lemmas $4.2-4.6$, we obtain
 \begin{equation}\left.\begin{split}
 G_{\geq 2t+3}(z,w)&=2\Re\left\{\left(\displaystyle\sum_{k=1}^{N}\Delta_{k}(z)\left(a_{k1}z_{1}+\dots+a_{kN}z_{N}\right)-(a+\overline{a})\Delta(z)\right)w^{t}\right\}
+2\Re\left\{P(z,\overline{z})\left(w^{t}-\langle
z,z\rangle^{t}\right)\right\}\\&\quad +2\Re\left\langle
\overline{z}+
\Delta'(z)+\Theta_{s}^{2}(z,\overline{z}),\overline{F_{\geq
2t+2}(z,w)}\right\rangle+\varphi'_{\geq2t+3}(z,\overline{z})-\varphi_{\geq2t+3}(z,\overline{z})
+\Theta_{ts+s+1}^{2t+3}(z,\overline{z}),
\end{split}\right.\label{51}\end{equation}
where $w$ satisfies (\ref{ec4}) and
$\mbox{wt}\left\{\overline{\Theta_{ts+s+1}^{2t+3}(z,\overline{z})}\right\}\geq
 ts+s+1$. It remains to study the
expression
\begin{equation}E(z,\overline{z})=2\Re\left\{P(z,\overline{z})\left(w^{t}-\langle
z,z \rangle^{t}\right)\right\}\label{E}.\end{equation}
 \bl For $w$
satisfying (\ref{ec4}) we make the following estimate
\begin{equation}E(z,\overline{z})=2\Re\left\{\left(P(z,\overline{z})+\overline{P(z,\overline{z})}\right)\Delta(z)
\displaystyle\sum_{k+l=t-1}w^{k}\langle
z,z\rangle^{l}\right\}+\Theta_{ts+s+1}^{2t+3}(z,\overline{z}),\end{equation}
where
$P(z,\overline{z})=\displaystyle\sum_{k,j=1}^{N}b_{k,j}z_{k}\overline{z}_{j}$
and
$\rm{wt}\left\{\overline{\Theta_{ts+s+1}^{2t+3}(z,\overline{z})}\right\}\geq
 ts+s+1$.
  \el
\begin{proof} We compute
\begin{equation}\left.\begin{split}E(z,\overline{z})&=2\Re\left\{P(z,\overline{z})\left(\Delta(z)+\overline{\Delta(z)}\right)\displaystyle\sum_{k+l=t-1}w^{k}\langle z,z\rangle^{l}\right\}
+\Theta_{ts+s+1}^{2t+3}(z,\overline{z}),\\&=2\Re\left\{\left(P(z,\overline{z})+\overline{P(z,\overline{z})}\right)\Delta(z)\displaystyle\sum_{k+l=t-1}w^{k}\langle
z,z\rangle^{l}\right\}
+\Theta_{ts+s+1}^{2t+3}(z,\overline{z}),\end{split}\right.
\end{equation}
where
$\mbox{wt}\left\{\overline{\Theta_{ts+s+1}^{2t+3}(z,\overline{z})}\right\}\geq
 ts+s+1$.
\end{proof}
We consider the following notations
\begin{equation}\left.\begin{split}& \quad\quad\quad\quad\quad\quad\quad\quad\mathcal{L}(z,\overline{z})=P(z,\overline{z})+\overline{P(z,\overline{z})}=\displaystyle\sum_{k,j=1}^{N}\left(b_{k,j}+\overline{b}_{j,k}\right)z_{k}\overline{z}_{j}
,\\&Q(z)=\displaystyle\sum_{k=1}^{N}\Delta_{k}(z)\left(a_{k1}z_{1}+\dots+a_{kN}z_{N}\right)-\left(a+\overline{a}\right)\Delta(z),
\quad
Q_{1}(z)=\displaystyle\sum_{k,j=1}^{N}\left(b_{k,j}+\overline{b}_{j,k}\right)z_{k}\Delta_{k}(z).\\&\end{split}\right.\label{KKK}\end{equation}
Then, for $w$ satisfying (\ref{ec4}), by Lemma $4.7$ and the
notations (\ref{KKK}), we can rewrite (\ref{51}) as follows
\begin{equation}\left.\begin{split}G_{\geq 2t+3}(z,w)& =2\Re\left\{Q(z)w^{t}\right\}
+2\Re\left\{\mathcal{L}(z,\overline{z})\Delta(z)E_{t-1}\left(w,\langle
z,z\rangle\right)\right\}+2\Re\left\langle \overline{z}+
\Delta'(z) +\Theta_{s}^{2}(z,\overline{z}),\overline{F_{\geq
2t+2}(z,w)}\right\rangle\\&\quad+\varphi'_{\geq
2t+3}(z,\overline{z})-\varphi_{\geq 2t+3}(z,\overline{z})
+\Theta_{ts+s+1}^{2t+3}(z,\overline{z}),
\end{split}\right.\label{53}\end{equation}
where
$\mbox{wt}\left\{\overline{\Theta_{ts+s+1}^{2t+3}(z,\overline{z})}\right\}\geq
 ts+s+1$. Here $E_{t-1}\left(w,\langle z,z
\rangle\right)=\displaystyle\sum_{k+l=t-1}w^{k}\langle
z,z\rangle^{l}$. For $p\geq 2t+3$ we prove the following lemma (
the analogue of  Lemma $3.4$ from Huang-Yin's paper ~\cite{HY2}):

 \bl We define $\epsilon\left(p\right)=0$ if $p<2t+s$ and $\epsilon\left(p\right)=1$  if $p\geq 2t+s$, $\gamma\left(p\right)=1$ if $p<ts+2$ and $\gamma(p)=0$ if
 $p=ts+2$. Let $N'_{s}:=ts+s+1$. For all $0\leq j\leq t$ and $p\in\left[2t+j(s-2)+3,2t+(j+1)(s-2)+2\right]$, we have the following estimate
\begin{equation}\left. \begin{split} G_{\geq p}(z,w)&=2(1-s)^{j}\Re\left\{Q(z)\Delta(z)^{j}w^{t-j}\right\}
+2\gamma\left(p\right)(-1)^{j}\Re\left\{\mathcal{L}(z,\overline{z})\Delta(z)^{j+1}\displaystyle\sum_{l_{1}+l_{2}=t-j-1}E_{l_{1},l_{2}}^{t-j}w^{l_{1}}\langle
z,z\rangle^{l_{2}}\right\}
\\&\quad+2\epsilon\left(p\right)\Re\left\{Q_{1}(z)\Delta(z)^{j}w^{t-j}\displaystyle\sum_{l=0}^{j-1}(-1)^{\beta_{l}}(1-s)^{l}F^{t-j}_{l}\right\}
+2\Re\left\langle \overline{z}+\Delta'(z)
+\Theta_{s}^{2}(z,\overline{z}),\overline{F_{\geq
p-1}(z,w)}\right\rangle \\&\quad+\varphi'_{\geq
p}(z,\overline{z})-\varphi_{\geq
p}(z,\overline{z})+\Theta_{N_{s}'}^{p}(z,\overline{z}),\end{split}\right.\label{54}\end{equation}
where
$\rm{wt}\left\{\overline{\Theta_{N'_{s}}^{p}(z,\overline{z})}\right\}\geq
 N'_{s}$ and $w$ satisfies (\ref{ec4}). Here $E_{l_{1},l_{2}}^{t-j}$ with $l_{1}+l_{2}=t-j-1$ and $F_{p}^{t-j}$ with $l=1,\dots,j-1$ are natural numbers
satisfying the following recurrence relations
$$F^{t-j-1}_{l+1}=F_{l}^{t-j},\quad F^{t-j-1}_{0}=\displaystyle\sum_{l_{1}+l_{2}
=t-j-1}E_{l_{1},l_{2}}^{t-j},\quad
E^{t-j-1}_{l,t-j-1-l}=\displaystyle\sum_{l'=l}^{t-j-l}E_{t-j-l',l'}^{t-j}.$$

 Also  $\beta_{l}\in\mathbb{N}$, for all $l=1,\dots,j-1$. \el
\begin{proof}

 For $j=0$ and $k=0$ we obtain $p=2t+3$. Therefore (\ref{54}) becomes
(\ref{53}).

 $\bf{Step\hspace{0.1 cm} 1.}$ We make a similarly approach as in the
Step $1$ of  Lemma $3.7$.

 $\bf{Step\hspace{0.1 cm} 2.}$ Assume that
we proved  Lemma $4.8$ for
$m\in\left[2t+j(s-2)+3,2t+(j+1)(s-2)+2\right]$, for $j\in[0,t-1]$.
We want to prove that (\ref{54}) holds for
$m\in\left[2t+(j+1)(s-2)+3,2t+(j+2)(s-2)+2\right]$. Collecting
from (\ref{54}) the terms of bidegree $(m,n)$ in
$(z,\overline{z})$ with $m+n=\Lambda+1:=2t+(j+1)(s-2)+2$, we
obtain
\begin{equation}\left. \begin{split} G_{nor}^{(\Lambda+1)}\left(z,\langle z,z\rangle\right)&
=2\Re\left\langle z,F_{nor}^{(\Lambda)}\left(z,\langle z,
z\rangle\right)\right\rangle+2\gamma\left(p\right)(-1)^{j}\Re\left\{\mathcal{L}(z,\overline{z})
\Delta(z)^{j+1}\langle
z,z\rangle^{t-j-1}\displaystyle\sum_{l_{1}+l_{2}=t-j-1}E_{l_{1},l_{2}}^{t-j}\right\}\\&\quad+
2\epsilon\left(p\right)\Re\left\{Q_{1}(z)\Delta(z)^{j}\langle z,z
\rangle^{t-j}\displaystyle\sum_{l=0}^{j-1}(-1)^{\beta_{l}}(1-s)^{l}F^{t-j}_{t-j}\right\}+2(1-s)^{j}\Re\left\{Q(z)\Delta(z)^{j}\langle
z,z\rangle^{t-j}\right\}
\\&\quad+\varphi'_{\Lambda+1}(z,\overline{z})
-\varphi_{\Lambda+1}(z,\overline{z})+\mathbb{P}_{N'_{s}}^{\Lambda+1}(z,\overline{z}),\end{split}\right.\label{55}\end{equation}
where
$\mbox{wt}\left\{\mathbb{P}^{\Lambda+1}_{N'_{s}}(z,\overline{z})\right\}\geq
 N'_{s}$. We define the following mappings
\begin{equation}\left.\begin{split}&\quad \quad\quad\quad\quad\quad\quad \quad\quad F_{1}^{(\Lambda)}(z,w)=-(1-s)^{j}Q(z)\Delta(z)^{j}w^{t-j-1}\left(z_{1},\dots,z_{N}\right),
\\&\quad\quad\quad\quad\quad F_{2}^{(\Lambda)}(z,w)=-\epsilon(p) Q_{1}(z)\Delta(z)^{j}w^{t-j-1}\left(\displaystyle\sum_{l=0}^{j-1}(-1)^{\beta_{l}}(1-s)^{l}F_{l}^{t-j}\right)\left(z_{1},\dots,z_{N}\right),\\&
F_{3}^{(\Lambda)}(z,w)=-\gamma(p)(-1)^{j}\Delta(z)^{j+1}\displaystyle\sum_{l_{1}+l_{2}=t-j-1}E_{l_{1},l_{2}}^{t-j}w^{t-j-1}\left(\displaystyle\sum_{l=1}^{N}\left(b_{l,1}+\overline{b}_{1,l}\right)z_{l},\dots,\displaystyle\sum_{l=1}^{N}
\left(b_{l,N}+\overline{b}_{N,l}\right)z_{l}\right),\\&
\quad\quad\quad\quad\quad\quad\quad\quad
F_{nor}^{(\Lambda)}(z,w)=F_{1}^{(\Lambda)}(z,w)+F_{2}^{(\Lambda)}(z,w)+F_{3}^{(\Lambda)}(z,w)+F_{4}^{(\Lambda)}(z,w),\end{split}\right.\label{56}\end{equation}
where $F_{4}^{(\Lambda)}(z,w)$ will be determined later (see
\ref{57}).

Substituting (\ref{56}) into (\ref{55}), by making some
simplifications it follows that
\begin{equation}\left. \begin{split}G_{nor}^{(\Lambda+1)}(z,\langle z,z\rangle)=&2\Re\left\langle z,F_{4}^{(\Lambda)}\left(z,\langle z,z\rangle\right)\right\rangle
+\varphi'_{\Lambda+1}(z,\overline{z})-\varphi_{\Lambda+1}
(z,\overline{z})+\mathbb{P}_{N'_{s}}^{\Lambda+1}(z,\overline{z}).\end{split}\right.\label{57}\end{equation}
By applying Extended Moser Lemma we find a solution
$\left(G_{nor}^{(\Lambda+1)}(z,w),F_{4}^{(\Lambda)}(z,w)\right)$
for (\ref{57}). By repeating the procedure from the first case of
the normal form construction, we obtain the following estimates
\begin{equation}\left.\begin{split}&\quad\quad\mbox{wt}\left\{G_{nor}^{(\Lambda+1)}(z,w)-G_{nor}^{(\Lambda+1)}\left(z,\langle
z,z\rangle\right)\right\} ,\hspace{0.1
cm}\mbox{wt}\left\{G_{nor}^{(\Lambda+1)}(z,w)\right\},\hspace{0.1
cm}\mbox{wt}\left\{ G_{nor}^{(\Lambda+1)}\left(z,\langle
z,z\rangle\right) \right\} \geq N'_{s},\\&\quad\quad\hspace{0.1
cm}\mbox{wt}\left\{F_{4}^{(\Lambda)}(z,w)-F_{4}^{(\Lambda)}\left(z,\langle
z,z\rangle\right)\right\},\hspace{0.1
cm}\mbox{wt}\left\{F_{4}^{(\Lambda)}(z,w)\right\},\hspace{0.1 cm}
\mbox{wt}\left\{F_{4}^{(\Lambda)}\left(z,\langle
z,z\rangle\right)\right\} \geq N'_{s}-s+1,\\&
\quad\quad\quad\quad\mbox{wt}\left\{
\overline{F_{4}^{(\Lambda)}(z,w)}\right\},\hspace{0.1
cm}\mbox{wt}\left\{\overline{F_{4}^{(\Lambda)}\left(z,\langle
z,z\rangle\right)}\right\},\hspace{0.1
cm}\mbox{wt}\left\{\overline{
F_{4}^{(\Lambda)}(z,w)-F_{4}^{(\Lambda)}\left(z,\overline{z}\right)}\right\}\geq
N'_{s}-1,\end{split}\right.\label{3112}\end{equation} where $w$
satisfies (\ref{ec4}). As a consequence of (\ref{3112}) we obtain
\begin{equation}\left.\begin{split}&\left\langle
\Delta'(z)+\Theta_{s}^{2}(z,\overline{z}),\overline{F_{4}^{(\Lambda)}(z,w)}\right\rangle+\left\langle
\overline{F_{4}^{(\Lambda)}(z,w)},\Delta'(z)+\Theta_{s}^{2}(z,\overline{z})\right\rangle
 =\Theta_{N'_{s}}^{\Lambda+2}(z,\overline{z})',\\&\quad\quad\quad
 \quad\quad
 \Re\left\langle
F_{4}^{(\Lambda)}(z,w)-F_{4}^{(\Lambda)}\left(z,\langle
z,z\rangle\right),z\right\rangle=\Theta_{N'_{s}}^{\Lambda+2}(z,\overline{z})',
\end{split}\right.\label{58}\end{equation}
where $w$ satisfies (\ref{ec4}) and each of
$\Theta_{N'_{s}}^{\Lambda+2}(z,\overline{z})'$ has the property
 $\mbox{wt}\left\{\overline{\Theta_{N'_{s}}^{2t+3}(z,\overline{z})}\right\}\geq
 N'_{s}$. Substituting
$F_{\geq\Lambda}(z,w)=F_{nor}^{(\Lambda)}(z,w)+F_{\geq
\Lambda+1}(z,w)$ and
$G_{\geq\Lambda+1}(z,w)=G_{nor}^{(\Lambda+1)}(z,w)+G_{\geq\Lambda+2}(z,w)$
in  (\ref{54}), it follows that
\begin{equation}\left.\begin{split} G_{nor}^{(\Lambda+1)}(z,w)+G_{\geq\Lambda+2}(z,w)=&2\Re\left\langle \overline{z}+ \Delta'(z)
+\Theta_{s}^{2}(z,\overline{z}),\overline{F_{nor}^{(\Lambda)}(z,w)}+\overline{F_{\geq\Lambda+1}(z,w)}\right\rangle
+\varphi'_{\Lambda+1}(z,\overline{z})-\varphi_{\Lambda+1}(z,\overline{z})\\&+\left(\Theta_{1}\right)_{N'_{s}}^{\Lambda+1}(z,\overline{z})
+\varphi'_{>\Lambda+1}(z,\overline{z})-\varphi_{>\Lambda+1}(z,\overline{z})
+\Theta_{N'_{s}}^{\Lambda+2}(z,\overline{z})\\&+2(1-s)^{j}\Re\left\{Q(z)\Delta(z)^{j}w^{t-j}\right\}\\&+
2\gamma(p)\Re\left\{(-1)^{j}\mathcal{L}(z,\overline{z})\Delta(z)^{j+1}\displaystyle\sum_{l_{1}+l_{2}=t-j-1}E_{l_{1},l_{2}}^{t-j}w^{l_{1}}\langle
z,z\rangle^{l_{2}}\right\}\\&+2\epsilon(p)\Re\left\{Q_{1}(z)\Delta(z)^{j}w^{t-j}\displaystyle
\sum_{l=0}^{j-1}(-1)^{\beta_{l}}(1-s)^{l}F^{t-j}_{l}\right\}
,\end{split}\right.\label{59}\end{equation} where $w$ satisfies
(\ref{ec4}). After a simplification in the preceding equation by
using (\ref{55}), it follows that
\begin{equation}G_{\geq\Lambda+2}(z,w)=2\Re\left\langle
\overline{z}+\Delta'(z)+\Theta_{s}^{2}(z,\overline{z}),\overline{F_{\geq\Lambda+1}(z,w)}\right\rangle
+\varphi_{\geq\Lambda+2}(z,\overline{z})-\varphi'_{\geq\Lambda+2}(z,\overline{z})+\Theta_{N'_{s}}^{\Lambda+2}(z,\overline{z})+J(z,\overline{z}),\label{60}\end{equation}
where we have used the following notation
\begin{equation}\left.\begin{split}J(z,\overline{z})&=2\Re\left\langle
z,F_{nor}^{(\Lambda)}(z,w)-F_{nor}^{(\Lambda)}\left(z,\langle
z,z\rangle\right)\right\rangle+2\Re\left\langle
\Delta'(z)+\Theta_{s}^{2}(z,\overline{z}),
\overline{F_{nor}^{(\Lambda)}(z,w)}\right\rangle\\&\quad
+2(1-s)^{j}\Re\left\{Q(z)\Delta(z)^{j}w^{t-j}-Q(z)\Delta(z)^{j}\langle
z,z\rangle^{t-j}\right\}+G_{nor}^{(\Lambda+2)}\left(z,\langle
z,z\rangle\right)-G_{nor}^{(\Lambda+2)}(z,w)\\&\quad+2\gamma(p)(-1)^{j}\Re\left\{\mathcal{L}(z,\overline{z})\Delta(z)^{j+1}\left(\displaystyle\sum_{l_{1}+l_{2}
=t-j-1}E_{l_{1},l_{2}}^{t-j}w^{l_{1}}\langle z,z\rangle^{l_{2}}
-\langle
z,z\rangle^{t-j-1}\displaystyle\sum_{l_{1}+l_{2}=t-j-1}E_{l_{1},l_{2}}^{t-j}\right)\right\}\\&\quad+2\epsilon(p)\Re\left\{Q_{1}(z)\Delta(z)^{j}\displaystyle
\sum_{l=0}^{j-1}(-1)^{\beta_{l}}(1-s)^{l}\left(F^{t-j}_{l}w^{t-j}-F^{t-j}_{l}\langle
z,z\rangle^{t-j} \right)\right\},\end{split}\right.\end{equation}
\begin{equation}\left.\begin{split}J(z,\overline{z})&=2\Re\left\langle
z,\displaystyle\sum_{k=1}^{3}\left(F_{k}^{(\Lambda)}(z,w)-F_{k}^{(\Lambda)}\left(z,\langle
z,z\rangle\right)\right)\right\rangle+2\Re\left\langle
\Delta'(z)+\Theta_{s}^{2}(z,\overline{z}),\displaystyle\sum_{k=1}^{3}\overline{F_{k}^{\left(\Lambda\right)}(z,w)}\right\rangle
\\&\quad+2(1-s)^{j}\Re\left\{Q(z)\Delta(z)^{j}\left(w^{t-j}-\langle
z,z\rangle^{t-j}\right)\right\}+G_{nor}^{(\Lambda+2)}\left(z,\langle
z,z\rangle\right)-G_{nor}^{(\Lambda+2)}(z,w)\\&\quad+2\gamma(p)(-1)^{j}\Re\left\{\mathcal{L}(z,\overline{z})\Delta(z)^{j+1}
\left(\displaystyle\sum_{l_{1}+l_{2}=t-j-1}E_{l_{1},l_{2}}^{t-j}w^{l_{1}}\langle
z,z\rangle^{l_{2}}-\displaystyle\sum_{l_{1}+l_{2}=t-j-1}E_{l_{1},l_{2}}^{t-j}\langle
z,z\rangle^{t-j-1}\right)\right\}\\&\quad+2\epsilon(p)\Re\left\{Q_{1}(z)\Delta(z)^{j}\displaystyle
\sum_{l=0}^{j-1}(-1)^{\beta_{l}}(1-s)^{l}F^{t-j}_{l}\left(w^{t-j}-\langle
z,z\rangle^{t-j}\right)\right\}.
\end{split}\right.\label{61}\end{equation}
We observe that
\begin{equation}\left.\begin{split}&\quad\quad\quad\quad
\Re \left\langle F_{1}^{(\Lambda)}\left(z,\langle
z,z\rangle\right),z
\right\rangle=-(1-s)^{j}\Re\left\{Q(z)\Delta(z)^{j}\langle
z,z\rangle^{t-j}\right\},\\&\quad\Re \left\langle
F_{2}^{(\Lambda)}\left(z,\langle z,z\rangle\right),z
\right\rangle=-\epsilon(p)\Re\left\{ Q_{1}(z)\Delta(z)^{j}\langle
z,z\rangle^{t-j}\displaystyle
\sum_{l=0}^{j-1}(-1)^{\beta_{l}}(1-s)^{l}F^{t-j}_{l}
\right\},\\&\hspace{0.1 cm}\Re \left\langle
F_{3}^{(\Lambda)}\left(z,\langle z,z\rangle\right),z
\right\rangle=-(-1)^{j}\gamma(p)\Re\left\{\mathcal{L}(z,\overline{z})\Delta(z)^{j+1}\langle
z,z\rangle^{t-j-1}\displaystyle\sum_{l_{1}+l_{2}=t-j-1}E_{l_{1},l_{2}}^{t-j}\right\}
\label{950}.\end{split}\right.\end{equation} Since
$\mbox{wt}\left\{F_{k}^{(\Lambda)}(z,w)\right\} \geq ts+1$ and
$\mbox{wt}\left\{\overline{F_{k}^{(\Lambda)}(z,w)}\right\} \geq
ts+s-1$ for all $k\in\{1,2,3\}$, it follows that
\begin{equation}2\Re\left\langle\Theta_{s}^{2}(z,\overline{z}),\displaystyle\sum_{k=1}^{3}\overline{F_{k}^{(\Lambda)}(z,w)}\right\rangle
=\Theta_{N'_{s}}^{\Lambda+2}(z,\overline{z}), \label{105}
\end{equation}
where
$\mbox{wt}\left\{\overline{\Theta_{N'_{s}}^{\Lambda+2}(z,\overline{z})}\right\}\geq
 N'_{s}$. By using (\ref{3112}), (\ref{58}), (\ref{950}), (\ref{105})  we can rewrite (\ref{61}) as follows
\begin{equation}\left.\begin{split}J(z,\overline{z})=&2\Re\left\langle
z,\displaystyle\sum_{k=1}^{3}F_{k}^{(\Lambda)}(z,w)\right\rangle+2\Re\left\langle
\Delta'(z),\displaystyle\sum_{k=1}^{3}\overline{F_{k}^{(\Lambda)}(z,w)}\right\rangle
\\&+2(1-s)^{j}\Re\left\{Q(z)\Delta(z)^{j}w^{t-j}\right\}+2(-1)^{j}\gamma(p)\Re\left\{\mathcal{L}(z,\overline{z})\Delta(z)^{j+1}\langle
z,z\rangle^{t-j-1}
\displaystyle\sum_{l_{1}+l_{2}=t-j-1}E_{l_{1},l_{2}}^{t-j}\right\}\\&+2\epsilon\left(p\right)\Re\left\{
Q_{1}(z)\Delta(z)^{j}w^{t-j}\displaystyle\sum_{l=0}^{j-1}(-1)^{\beta_{l}}(1-s)^{l}F^{t-j}_{l}\right\}.
\end{split}\right.\label{700}\end{equation}
Substituting the formulas  of $F_{1}^{(\Lambda)}(z,w)$,
$F_{2}^{(\Lambda)}(z,w)$ and $F_{3}^{(\Lambda)}(z,w)$ in
(\ref{700}) and  using $w$ satisfying (\ref{ec4}), we obtain
\begin{equation}\left.\begin{split}J(z,\overline{z})&=-2(1-s)^{j}\Re\left\{Q(z)\Delta(z)^{j}w^{t-j-1}\left(\left\langle
z,z\right\rangle+s\Delta(z)\right)
-Q(z)\Delta(z)^{j}w^{t-j}\right\}\\&\quad-2(-1)^{j}\gamma(p)\Re\left\{\mathcal{L}(z,\overline{z})\Delta(z)^{j+1}
\displaystyle\sum_{l_{1}+l_{2}
=t-j-1}E_{l_{1},l_{2}}^{t-j}w^{l_{2}}\left(w^{l_{1}}-\langle
z,z\rangle^{l_{1}}\right)\right\}
\\&\quad-2(-1)^{j}\gamma(p)\Re\left\{Q_{1}(z)\Delta(z)^{j+1} \displaystyle\sum_{l_{1}+l_{2}
=t-j-1}E_{l_{1},l_{2}}^{t-j}w^{t-j-1}\right\}\\&\quad-2\epsilon(p)\Re\left\{Q_{1}(z)\Delta(z)^{j}\displaystyle\sum_{l=0}^{j-1}(-1)^{\beta_{l}}(1-s)^{l}
F_{l}^{t-j}w^{t-j-1}\left(\langle z,z\rangle+s\Delta(z)-w\right)
\right\},\end{split}\right.\label{7709}\end{equation}

 By (\ref{7709}) and by the next identity (\ref{709}) we obtain the recurrence
relations given by the statement of Lemma $23$.

\begin{equation}\left.\begin{split}J(z,\overline{z})&=2(1-s)^{j+1}
\Re\left\{Q(z)\Delta(z)^{j+1}w^{t-j-1}\right\}\\&\quad+2\gamma(p)(-1)^{j+1}\Re\left\{\mathcal{L}(z,\overline{z})\Delta(z)^{j+2}
\displaystyle\sum_{l_{1}+l_{2}
=t-j-2}E_{l_{1},l_{2}}^{t-j-1}w^{l_{2}}\langle z,z
\rangle^{l_{1}}\right\}\\&\quad+2(-1)^{j+1}\Re\left\{Q_{1}(z)\Delta(z)^{j+1}
\displaystyle\sum_{l_{1}+l_{2}
=t-j-1}E_{l_{1},l_{2}}^{t-j}w^{t-j-1}\right\}
\\&\quad+2\epsilon(p)\Re\left\{Q_{1}(z)\Delta(z)^{j+1}\displaystyle\sum_{l=0}^{j-1}(-1)^{\beta_{l}+1}(1-s)^{l+1}F_{l}^{t-j}w^{t-j-1}\right\}+\Theta_{N'_{s}}^{\Lambda+2}(z,\overline{z}),
\end{split}\right.\label{709}\end{equation}
where
$\mbox{wt}\left\{\overline{\Theta_{N'_{s}}^{\Lambda+2}(z,\overline{z})}\right\}\geq
 N'_{s}$.

 The proof of our Lemma follows by using (\ref{709}) and (\ref{60}).
 \end{proof}

Collecting the terms of bidegree $(m,n)$ in $(z,\overline{z})$
from (\ref{54}) with $m+n=ts+s$ and $t=j$, we obtain
\begin{equation}\left.\begin{split}
 G_{nor}^{(ts+s)}(z,\langle z,z\rangle)=&2(1-s)^{t}\Re\left\{Q(z)\Delta(z)^{t}\right\}+2K\Re\left\{Q_{1}(z)\Delta(z)^{t}\right\}+2\Re\left\langle
 z,F_{nor}^{(ts+s-1)}(z,w)\right\rangle\\&+\varphi'_{ts+s,0}(z,\overline{z})-\varphi_{ts+s,0}(z,\overline{z})
 +\left(\Theta_{1}\right)_{N'_{s}}^{ts+s}(z,\overline{z}).\end{split}\right.\label{65}\end{equation}
By applying Extended Moser Lemma we find a solution
$\left(G_{nor}^{(ts+s)}(z,w),F_{nor}^{(ts+s-1)}(z,w)\right)$ for
(\ref{65}). Collecting the pure terms of degree $ts+s$ from
(\ref{65}), it follows that
\begin{equation}\varphi'_{ts+s,0}(z)-\varphi_{ts+s,0}(z)=(1-s)^{t}Q(z)\Delta(z)^{t}+KQ_{1}(z)\Delta(z)^{t},\label{VR}\end{equation}
where
 $K=(-1)^{\beta_{1}}k_{1}(1-s)^{t-1}+\dots+(-1)^{\beta_{t-1}}k_{t-1}(1-s)+(-1)^{\beta_{t}}k_{t}$,
with $k_{1},\dots ,k_{t}\in\mathbb{N}$. By the proof of  Lemma
$4.8$ (see (\ref{7709}) and (\ref{709})) we observe that
$\beta_{1}=1,\dots,\beta_{t}=t$. Next, by applying Lemma $2.4$ to
$\varphi_{ts+s,0}(z)$ and
 $\varphi'_{ts+s,0}(z)$, it follows that
\begin{equation}\left.\begin{split}&\hspace{0.1
cm}\varphi_{ts+s,0}(z)=\left(A_{1}(z)\Delta_{1}(z)+\dots+A_{N}(z)\Delta_{N}(z)\right)\Delta(z)^{t}+C(z),
\\&\hspace{0.1
cm}\varphi'_{ts+s,0}(z)=\left(A'_{1}(z)\Delta_{1}(z)+\dots+A'_{N}(z)\Delta_{N}(z)\right)\Delta(z)^{t}+C'(z),\end{split}\right.
\label{fin}\end{equation} where
$\left(\Delta_{k}\Delta^{t}\right)^{\star}\left(C(z)\right)=\left(\Delta_{k}\Delta^{t}\right)^{\star}\left(C'(z)\right)=0$,
for all $k=1,\dots,N$. We have
\begin{equation}\left.\begin{split}&\quad\quad\quad\quad Q(z)=\displaystyle\sum_{k=1}^{N}\Delta_{k}(z)\left(a_{k1}z_{1}+\dots+\left(a_{kk}-\frac{a+\overline{a}}{s}\right)z_{k}
+\dots+a_{kN}z_{N}\right),\\&Q_{1}(z)=\displaystyle\sum_{k=1}^{N}\Delta_{k}(z)\left(\left(a_{k1}+\overline{a}_{1k}\right)z_{1}+\dots+\left(a_{kk}+\overline{a}_{kk}-(a+\overline{a})\right)z_{k}
+\dots+\left(a_{kN}+\overline{a}_{Nk}\right)z_{N}\right).\end{split}\right.\end{equation}
We impose the normalization condition
$\left(\Delta_{k}\Delta^{t}\right)^{\star}\left(\varphi_{ts+s,0}'(z)\right)=0$,
for all $k=1,\dots,N$.
 By Lemma $2.4$ this is equivalent to finding  $\left(a_{ij}\right)_{1\leq i,j \leq
N}$ such that $A'_{1}(z)=\dots=A'_{N}(z)=0$.  It follows that
\begin{equation}\left.\begin{split}&\quad\quad(1-s)^{t}a_{kj}+K\left(a_{kj}+\overline{a}_{jk}\right)=c_{kj},\quad\mbox{for all}\hspace{0.1 cm}k,j=1,\dots, N,\hspace{0.1 cm}k\neq j,
\\&(1-s)^{t}\left(a_{kk}-\frac{a+\overline{a}}{s}\right)+K\left(a_{kk}+\overline{a}_{kk}-(a+\overline{a})\right)=c_{kk},\quad\mbox{for all}\hspace{0.1 cm}k=1,\dots, N,\end{split}\right.\label{470}\end{equation}
where $c_{kj}$ is determined, for all $k,j=1,\dots,N$. Here
$Na=\displaystyle\sum_{k=1}^{N}a_{kk}$. By using the second
equation from (\ref{470}) we find $\Im a_{kk}$, for all
$k=1,\dots,N$. By taking the real part in the second equation from
(\ref{470}), we obtain
\begin{equation}\left(Ns(1-s)^{t}+2NKs\right)\Re
a_{kk}-\left(2(1-s)^{t}+2Ks\right)\displaystyle\sum_{l=1}^{N}\Re
a_{ll}=\Re c_{k,k}, \quad k=1,\dots,N.\label{yu}\end{equation} By
summing all the identities from (\ref{yu}), it follows that
$(1-s)^{t}N(s-2)\displaystyle\sum_{l=1}^{N}\Re
a_{ll}=\displaystyle\sum_{k=1}^{N}\Re c_{k,k}$. Next, going back
to (\ref{yu}) we  find $\Re a_{ll}$, for all $l=1,\dots,N$. Now,
let $k\neq j$ and $k,j\in\{1,\dots, N\}$. By taking the real and
the imaginary part in first equation from (\ref{470}), we obtain
\begin{equation}\left.\begin{split}&\hspace{0.1 cm}\left((1-s)^{t}+K\right)\Re a_{kj}+K\Re a_{jk}=\Re c_{k,j},\quad K\Re a_{kj}+\left((1-s)^{t}+K\right)\Re a_{jk}=\Re c_{j,k},
\\&\left((1-s)^{t}+K\right)\Im a_{kj}-K\Im a_{jk}=\Im c_{k,j},\quad -K\Im a_{kj}+\left((1-s)^{t}+K\right)\Im a_{jk}=\Im c_{j,k},\label{473}\end{split}\right.\end{equation}
where $c_{k,j}$ is determined, for all $k,j=1,\dots,N$ and $k\neq
j$. In order to solve the preceding system of equations it is
enough to observe that $(1-s)^{t}\left((1-s)^{t}+2K\right)\neq 0$.
It is equivalent to observe that
\begin{equation}\left.\begin{split}&\quad\quad\hspace{0.2 cm}(1-s)^{t}+2\left((-1)k_{1}(1-s)^{t-1}+\dots+(-1)^{t}k_{t}\right)\neq 0,
\\&(-1)^{t}\left((s-1)^{t}+2\left(k_{1}(s-1)^{t-1}+k_{2}(s-1)^{t-2}+\dots+k_{t}\right)\right)\neq 0.\end{split}\right.\end{equation}

\vspace{0.1 cm}

 By composing the map that sends $M$ into
(\ref{ec4}) with the map (\ref{48}) we obtain our formal
transformation that sends $M$ into $M'$ up  to degree $ts+s+1$.
\section{Proof of Theorem $5$-Uniqueness of the formal transformation map}

 In order to prove the uniqueness of the map (\ref{2800}) it is enough to prove that the
following map is the identity
\begin{equation}M'\ni (z,w)\longrightarrow\left(z+\displaystyle\sum_{k\geq2}F_{nor}^{(k)}(z,w),w+\displaystyle\sum_{k
\geq2}G_{nor}^{(k+1)}(z,w)\right)\in M'.\label{XXX1}\end{equation}
Here $M'$ is a manifold defined by the normal form from the
Theorem $1.5$. We have used the notations (\ref{2.11}). We perform
induction on $k\geq 2$. \bd The undetermined homogeneous parts of
the map (\ref{XXX1}) by applying Extended Moser Lemma are called
the free parameters. \ed

We prove that $F_{nor}^{(2)}(z,w)=0$. Here we recall the first
case of the normal form construction. We assume that $t=1$. By
repeating the normalization procedures from the first case of the
normal form construction, we  find that all of the homogeneous
components of $F_{nor}^{(2)}(z,w)$ except the free parameter  are
$0$ and that $G_{nor}^{(3)}(z,w)=0$. By using the same approach as
in the first case of the normal form construction (see
(\ref{f01})), it follows that
\begin{equation}\varphi'_{s+1,0}(z)-\varphi_{s+1,0}(z)=(1-s)\langle z,a\rangle\Delta(z)=0.\label{u1}\end{equation}
Here $a$ is the free parameter of $F_{nor}^{(2)}(z,w)$. It follows
that $a=0$. Therefore $F_{nor}^{(2)}(z,w)=0$.

 We assume that $F_{nor}^{(2)}(z,w)=\dots=F_{nor}^{(k-2)}(z,w)=0$, $G_{nor}^{(3)}(z,w)=\dots=G_{nor}^{(k-1)}(z,w)=0$. We want to prove that
$F_{nor}^{(k-1)}(z,w)=0$, $G_{nor}^{(k)}(z,w)=0$. First, we
consider the case when $k=2t$, with $t\geq 2$. Let
$a\in\mathbb{C}^{N}$ be the free parameter of the polynomial
$F_{nor}^{(2t)}(z,w)$. By repeating all the normalization
procedures from the first case of the normal form construction  it
follows that  all of the homogeneous components of
$F_{nor}^{(2t)}(z,w)$
 except the free parameters are $0$ and that $G_{nor}^{(2t+1)}(z,w)=0$.
We are interested in the image of the manifold $M$ through  the
map (\ref{XXX1}) to $M$ up to degree $ts+1$. We repeat the
normalization procedure done during Lemma $3.7$ proof. In that
case we have considered
 a particular mapping (see (\ref{t1})). Here we have a general polynomial map with other free parameters. They generate terms
 of weight at least $ts+2$ that do not change their weight under the conjugation:
    \begin{equation}\left.\begin{split}&\mbox{wt}\left\{\left\langle F_{1,m}(z)w^{m},z \right \rangle\right\},\quad\mbox{wt}\left\{\left\langle z,F_{1,m}(z)w^{m} \right \rangle\right\}\geq ts+2,\hspace{0.1 cm}\mbox{for all}\hspace{0.1 cm}m>t;
     \\&\mbox{wt}\left\{\left\langle F_{0,r}(z) w^{r},z\right\rangle\right\},\quad\mbox{wt}\left\{\left\langle z,F_{0,r}(z) w^{r}\right\rangle\right\}\geq ts+2,\hspace{0.1 cm}\mbox{for all}\hspace{0.1 cm}r\geq
     t+2.\end{split}\right.\end{equation}
Here $F_{1,m}(z)w^{m}$, $F_{0,r}(z) w^{r}$ are the free parameters
of $F_{nor}^{(2m+1)}(z,w)$ and $F_{nor}^{(2r)}(z,w)$, for all
$m>t$ and $r\geq t+2$. Therefore they cannot interact with the
pure terms of degree $ts+1$ (because of the higher weight). All
the Lemmas $3.1-3.6$ remain the same in this general case.

 By using the same approach as in the first case of the normal form
construction (see (\ref{caz1})), it follows that
\begin{equation}\varphi'_{ts+1,0}(z)-\varphi_{ts+1,0}(z)=(1-s)^{t}\langle z,a\rangle\Delta^{t}(z)=0.\label{u1}\end{equation}
 It follows that $a=0$. Therefore $F_{nor}^{(2t)}(z,w)=0$.

We assume that $k=2t+1$, with $t\geq 2$. Let
$\left(a_{i,j}\right)_{1\leq i,j \leq N}$  be the free parameter
of $F_{nor}^{(2t+1)}(z,w)$. By repeating all the normalization
procedures from the first case of the normal form construction, it
follows that
 all of the homogeneous components of $F_{nor}^{(2t+1)}(z,w)$
 except the free parameters are $0$ and that $G_{nor}^{(2t+2)}(z,w)=0$.

 We are interested of the image of the manifold $M'$ through the map
(\ref{t1}) to $M'$ up to degree $ts+s+1$. The other free
parameters of the map (\ref{XXX1}) generate terms
 of weight at least $ts+s+1$ that do not change their weight under the conjugation:
      \begin{equation}\left.\begin{split}&\mbox{wt}\left\{\left\langle F_{1,m}(z)w^{m},z \right \rangle\right\},\quad\mbox{wt}\left\{\left\langle z,F_{1,m}(z)w^{m} \right \rangle\right\}\geq ts+s+1,\hspace{0.1 cm}\mbox{for all}\hspace{0.1 cm}m>t+1;\\&
    \quad \mbox{wt}\left\{\left\langle F_{0,r}(z) w^{r},z\right\rangle\right\},\quad\mbox{wt}\left\{\left\langle z,F_{0,r}(z) w^{r}\right\rangle\right\}\geq ts+s+1,\hspace{0.1 cm}\mbox{for all}\hspace{0.1 cm}r\geq
     t+3.\end{split}\right.\end{equation}
All  Lemmas $4.1-4.7$ remain true in this general case.

  By using the same approach as in the second case of the normal form construction (see (\ref{VR})), it follows that
\begin{equation}\varphi'_{ts+s,0}(z)-\varphi_{ts+s,0}(z)=(1-s)^{t}Q(z)\Delta^{t}(z)=0.\label{u1}\end{equation}
 It follows that
$\left(a_{i,j}\right)_{1\leq i,j \leq N}=0$. Therefore
$F_{nor}^{(2t+1)}(z,w)=0,\quad G_{nor}^{(2t+2)}(z,w)=0.$ This
proves that  (\ref{XXX1}) is the identity mapping. From here we
conclude the uniqueness of the formal transformation (\ref{2800}).









\begin{thebibliography}{BER96b}

\bibitem{AG} {\bf Ahern,~P.; Gong,~X.} --- Real analytic manifolds in $\mathbb{C}^{n}$ with parabolic complex tangents along
a submanifold of codimension one. {\em Ann. Fac. Sci. Toulouse
Math. (6)} {\bf 18} (2009), no. 1, 1--64.
\bibitem{B} {\bf Bishop,~E.} --- Differentiable Manifolds In Complex
Euclidian Space. {\em Duke Math. J.} {\bf 32} (1965), no. 1,
1--21.
\bibitem{CM} {\bf Chern,~S.S; Moser,~J.K.} --- Real hypersurfaces in complex manifolds. {\em Acta Math.} {\bf 133},
(1974), no. 1, 219--271.
\bibitem{Cof1} {\bf Coffman,~A.} --- Analytic Normal Form For CR Singular
Surfaces in $\mathbb{C}^{3}$. {\em Houston J. Math.} {\bf 30}
(2004), no. $4$, 969-996.
\bibitem{Cof2} {\bf Coffman,~A.}--- CR Singularities of
real threefolds in  $\mathbb{C}^{4}$. {\em Adv. Geom.} {\bf 6}
(2006), no 1, 109-137.
\bibitem{Cof3} {\bf Coffman,~A.}--- Unfolding CR Singularities
 {\em Mem. Amer. Math. Soc.} {\bf 205} (2010), no. 962.
\bibitem{Cof4} {\bf Coffman,~A.}--- CR singularities of
real fourfolds in $\mathbb{C}^{3}$. {\em Illinois J. Math.} {\bf
53} (2009), no. 3, 939-981.
\bibitem{DTZ} {\bf Dolbeault,~ P.; Tomassini,~G.; Zaitsev,~D. }--- On Levi-flat hypersurfaces with prescribed boundary.
 {\em Pure and Applied Mathematics Quarterly} {\bf 6} (2010), no 3, ({\em Special Issue: In honor of Joseph J. Kohn. Part 1}),
 725--753.
 \bibitem{DTZ1} {\bf Dolbeault,~ P.; Tomassini,~G.; Zaitsev,~D. }--- Boundary Problem for
  Levi-Flat Graphs. {\em Indiana Univ. Math. J.} {\bf 60} (2011), no. 1,
161--170.
\bibitem{E2} {\bf Ebenfelt,~P. }--- New invariants in CR structures and a normal form for real hypersurfaces at a Levi degeneracy.
{\em J. Differential Geom. }{\bf 50}, (1998), no 2, 207-247.
\bibitem{G1} {\bf Gong,~X.} --- Normal forms of real surfaces under unimodular transformations near elliptic complex tangents. {\em Duke. Math. J.} {\bf 74} (1994), no. 1, 145--157.
\bibitem{G2} {\bf Gong,~X.} --- On the convergence of normalizations of real analytic surfaces near hyperbolic complex tangents. {\em Comment. Math. Helv.} {\bf 69} (1994), no. 4, 549--574.
\bibitem{G3} {\bf Gong,~X.} --- Existence of real analytic surfaces with hyperbolic complex tangent that are formally,  but not holomorphically equivalent to quadrics. {\em Indiana Univ.
Math. J.} {\bf 53} (2004), no. 1, 83--95.
\bibitem{HY2} {\bf Huang,~X.; Yin,~W.} --- A Bishop surface with vanishing Bishop invariant. {\em Invent. Math.}
 {\bf 176} (2010), no 3, 461-520.
\bibitem{HY1} {\bf Huang,~X.; Yin,~W.} --- A codimension two CR singular submanifold that is formally equivalent
 to a symmetric quadric. {\em Int. Math. Res. Notices IMRN}
  (2009), no 15, 2789-2828.
\bibitem{AHbook} {\bf Huang,~X.} --- {\em Local Equivalence Problems for Real Submanifolds in Complex Spaces}.
 Lecture Notes in Mathematics, Springer-Verlag,
pp. 109-161, Berlin-Heidelberg-New York, (2004).
\bibitem{HY4} {\bf Huang,~X.; Yin,~W.} --- Equivalence problem for Bishop surfaces. {\em Sci. China Math.}
 {\bf 43} (2010), no 3, 687-700.
\bibitem{HY5} {\bf Huang,~X.} --- On a $n$-manifold in $\mathbb{C}^{n}$ near an elliptic complex tangent. {\em J. Amer. Math. Soc.}
 {\bf 11} (1998), no 3, 669-692.
  \bibitem{HY6} {\bf Huang,~X.; Krantz,~S.} --- On a problem of Moser. {\em Duke Math. J.}
 {\bf 78} (1995), no 1, 213-228.
\bibitem{Mos} {\bf  Moser,~J.} ---Analytic Surfaces in $\mathbb{C}^{2}$ and their local hull
of holomorphy. {\em Ann. Acad. Sci. Fenn. Ser. A.I. Math.} {\bf
10} (1985), 397-410.
\bibitem{MW} {\bf Moser,~ J.; Webster,~S.} --- Normal forms for real surfaces in
$\mathbb{C}^{2}$ near complex tangents and hyperbolic surface
transformations. {\em Acta Math.} {\bf 150} (1983), 255--296.
  \bibitem{L1} {\bf Lebl,~J.} --- Nowhere minimal CR submanifolds and Levi-flat hypesurfaces. {\em J. Geom. Anal.}
 {\bf 17} (2007), no 2, 321-341.
  \bibitem{L2} {\bf Lebl,~J.} --- Extension of Levi-flat hypersurfaces past CR boundaries. {\em Indiana Univ. Math. J.}
 {\bf 57} (2008), no 2, 699-716.
\bibitem{L3} {\bf Lebl,~J.} --- Levi-flat hypersurfaces with real analytic boundary. {\em Trans. Amer. Math. Soc.}
 {\bf 362} (2010), no 12, 6367-6380.
\bibitem{Sh} {\bf Shapiro,~H.} ---Algebraic Theorem
of E.Fischer and the holomorphic Goursat problem. {\em Bull. London
Math. Soc.} {\bf 21}  (1989), no 6, 513--537.
\bibitem{DZ1} {\bf Zaitsev,~D.} --- New Normal Forms for Levi-nondegenerate
Hypersurfaces. {\em Several Complex Variables and Connections with
PDE Theory and Geometry}. Complex analysis-Trends in Mathematics,
Birkhauser Verlag, ({\em Special Issue: In the honor of Linda
Preiss Rothschild}), pp. 321-340, Basel/ Switzerland, (2010). {\em
$http://arxiv.org/abs/0902.2687.$}
\bibitem{DZ2} {\bf Zaitsev,~D.} --- Normal forms of non-integrable almost CR structures. {\em Amer. J. Math.} {\bf
134} (2012),no.4, 915-947.

\end{thebibliography}
\end{document}